\documentclass{rsproca}

\usepackage{paralist}
\usepackage{subcaption}
\usepackage{epstopdf}
\usepackage[square,numbers]{natbib}

\graphicspath{{Figures/}}

\newcommand{\R}{\mathbb{R}}
\renewcommand{\d}{\mathrm{d}}

\newcommand{\y}{\mathbf{y}}

\newcommand{\w}{\mathbf{w}}
\newcommand{\Asys}{A_\mathrm{sys}}
\renewcommand{\v}{\mathbf{v}}

\definecolor{dgreen}{rgb}{0,0.5,0}
\definecolor{dblue}{rgb}{0,0,0.5}
\newcommand{\pr}[1]{{\color{black}#1}}
\newcommand{\del}[1]{\iffalse{\color{dgreen}#1}\fi} 
\newcommand{\js}[1]{{\color{black}#1}} 
 
\newcommand{\prnew}[1]{{\color{black}#1}} 

\newcommand{\comment}[1]{\iffalse{\color{dgreen}[#1]}\fi} 
\newcommand{\prrsa}[1]{{\color{black}#1}}

\newcommand{\changed}[1]{{\color{black}#1}}

\newcommand{\finalchange}[1]{{\color{black}#1}}

\jname{rspa}
\Journal{Proc R Soc A\ }

\begin{document}

\title{Inverse-square law between time and amplitude for crossing
  tipping thresholds}

\author{
Paul Ritchie$^{1}$, {\"O}zkan Karabacak$^{2}$ and Jan Sieber$^{3}$}

\address{$^{1}$Earth System Science, College of Life and Environmental Sciences, Harrison Building, University of Exeter, Exeter, EX4 4QF, United Kingdom\\
$^{2}$Department of Electronic Systems, Automation and Control, Aalborg University, Fredrik Bajers Vej 7 C, 9220 Aalborg East, Denmark\\
$^{3}$Centre for Systems, Dynamics and Control, College of Engineering, Mathematics and Physical Sciences, Harrison Building, University of Exeter, Exeter, EX4 4QF, United Kingdom}

\subject{Tipping Points}

\keywords{tipping point, overshoot, bifurcation}

\corres{Paul Ritchie\\
\email{Paul.Ritchie@exeter.ac.uk}}

\begin{abstract}
  A classical scenario for tipping is that a dynamical system
  experiences a slow parameter drift across a fold tipping point,
  caused by a run-away positive feedback loop. We study \prrsa{what
    happens if one turns around after} one has crossed the
  threshold. We derive a simple criterion that relates \prrsa{how far
    the parameter exceeds the tipping threshold maximally and how long
    the parameter stays above the threshold to avoid tipping} in an inverse-square law
  to observable properties of the dynamical system near the fold.

  For the case when the dynamical system is subject to stochastic
  forcing we give an approximation to the probability of tipping \prrsa{if
  a parameter changing in time reverses} near the tipping point.

    
The derived approximations are valid if the parameter \prrsa{change in time is}
sufficiently slow. We demonstrate for a higher dimensional
system, a model for the Indian summer monsoon, how numerically
observed escape from the equilibrium converge to our asymptotic
expressions.  The inverse-square law between peak of the parameter
\prrsa{forcing and the time the parameter spends above a given
  threshold} is also visible in the level curves of equal probability when
the system is subject to random disturbances.
\end{abstract}


\begin{fmtext}

\end{fmtext}


\maketitle

\section{Introduction}
\label{sec:Intro}

The phenomenon of tipping is subject to ongoing intense study within
the scientific community due to its prominence in complex systems,
including climate
\citep{lenton2008tipping,held2004detection,holland2006future,boulton2014early},
ecosystems \citep{laurance201110,clark2013light,siteur2014beyond,gandhi2015localized} and finance
\citep{yan2010diagnosis}. The notion of tipping usually refers to a
sudden large qualitative change in output behavior caused by a small
change to input levels or rates \citep{ashwin2017parameter}. The
classical and most common model case for tipping is that the system
can be described (possibly at a coarse level) as a dynamical system
with a slowly drifting system parameter which passes slowly through a
fold (or saddle-node) bifurcation. In scientific terms the
mathematical scenario of a fold bifurcation at some system parameter
value corresponds to the presence of internal positive feedback loops,
which, with sufficient internal amplification, lead to a run-away
scenario. 

In Section \ref{sec:Monsoon} we will introduce a simple \finalchange{conceptual} Indian summer
monsoon model, originally derived by \citet{zickfeld2004modeling}. In this model a
positive feedback loop is formed between the temperature difference
over the Indian Ocean and Indian subcontinent and moisture advection
\citep{levermann2009basic}. In the summer months the temperature over
land warms quicker than the temperature over the ocean, which creates
winds coming off the ocean onto land \citep{zickfeld2004modeling}. The
winds carry moisture which is deposited over the land in the form of
precipitation. This process releases latent heat, causing the
temperature over land to increase, creating a greater temperature
difference and thus generating stronger winds to complete the positive
feedback loop. \citet{zickfeld2005indian} identified a tipping
threshold in the planetary albedo \prrsa{(the fraction of incoming
  solar radiation that is reflected over the Indian subcontinent)}, such that
increasing the albedo above this value will cause the monsoon to
shutdown.

The classical tipping scenario considers a gradual parameter change
that varies the system parameter slowly through the tipping threshold
(the fold bifurcation parameter value), causing a transition from the
current equilibrium, gradually varying with the parameter, to a new
state, possibly far away in state space. However, we may expect that
this transition is delayed with respect to the passage through the
tipping threshold if the system is forced at a faster than
infinitesimal speed
\citep{berglund2006noise,majumdar2013transitions}. This delay may pose
policy relevant questions, since many real life scenarios,
particularly in climate \citep{lenton2011early}, display similar
characteristics to those of a fold bifurcation transgression. For
example, \finalchange{numerical simulations of climate models suggest that} the Atlantic Meridional Overturning Circulation (AMOC) can be
disrupted or even stopped by an increase of freshwater to the North
Atlantic \citep{hawkins2011bistability}. However, due to the slow
response time of the system \citep{zhu2015amoc,li2016atlantic} it may
be possible to exceed the critical threshold for some time but still
maintain the circulation if the freshwater forcing is reduced to
values below a critical level sufficiently rapidly.

We will investigate this for the example model of the Indian summer
monsoon, one of the policy relevant tipping elements in the climate
system identified by \citet{lenton2008tipping}. For this climate
subsystem, policy makers may be interested in understanding: if the
albedo was increased beyond the threshold, can the albedo be reversed
quickly enough to prevent a shutdown of the monsoon? We use Zickfeld's
model to illustrate the \prrsa{deterministic} inverse square law for
the maximal permitted exceedance value and time over the tipping
threshold, \prrsa{and the deviations from it affected by random
disturbances.}

\prrsa{\paragraph*{Deterministic result} 
\begin{figure}[h!]
        \centering
        \subcaptionbox{\label{Albedo time profiles:intro}}[0.32\linewidth]
                {\includegraphics[scale = 0.3]{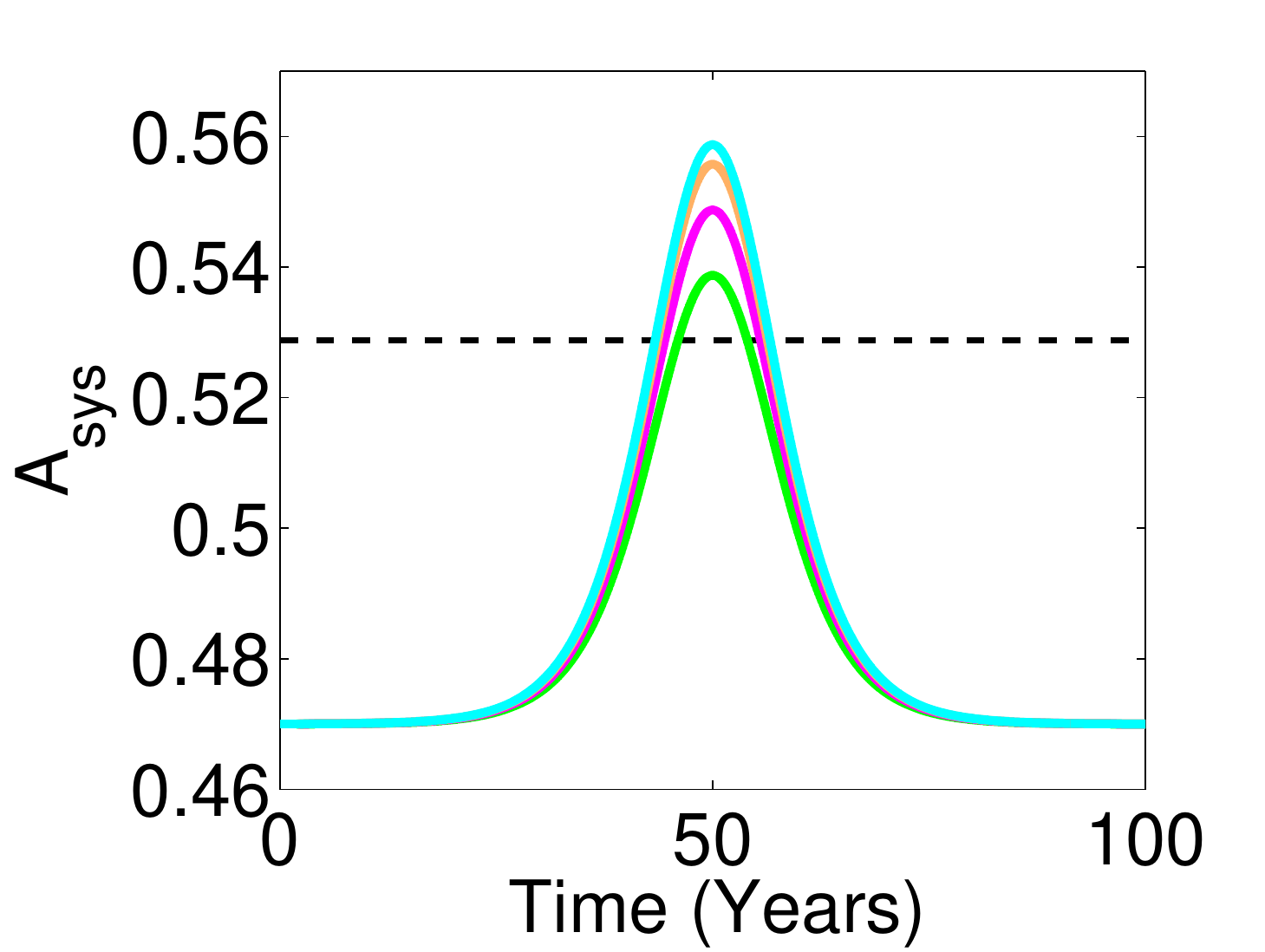}}
       \hfill
        \subcaptionbox{\label{Humidity time profiles:intro}}[0.32\linewidth]
                {\includegraphics[scale = 0.3]{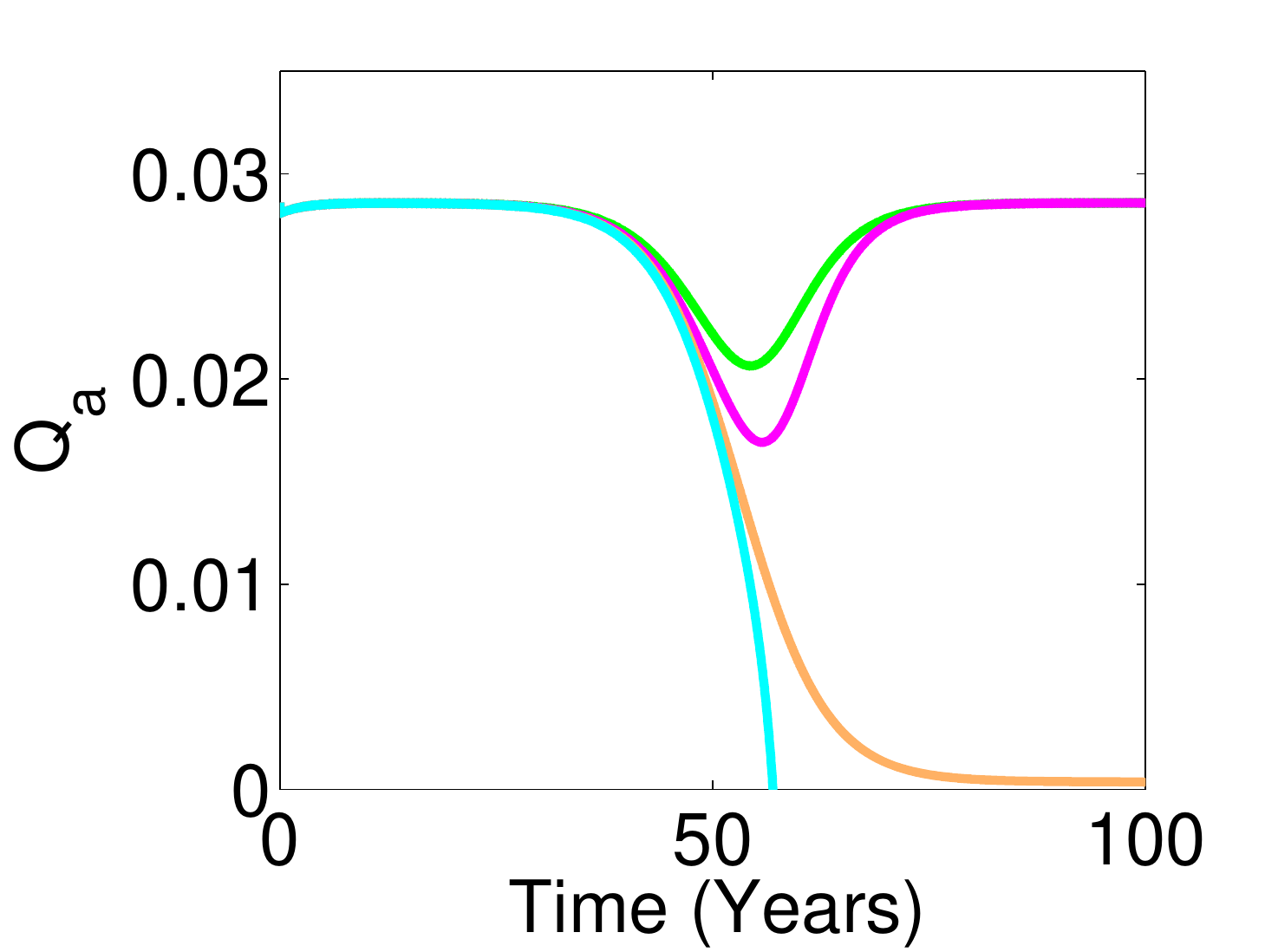}}
       \hfill
        \subcaptionbox{\label{Monsoon tajectories:intro}}[0.32\linewidth]
                {\includegraphics[scale = 0.3]{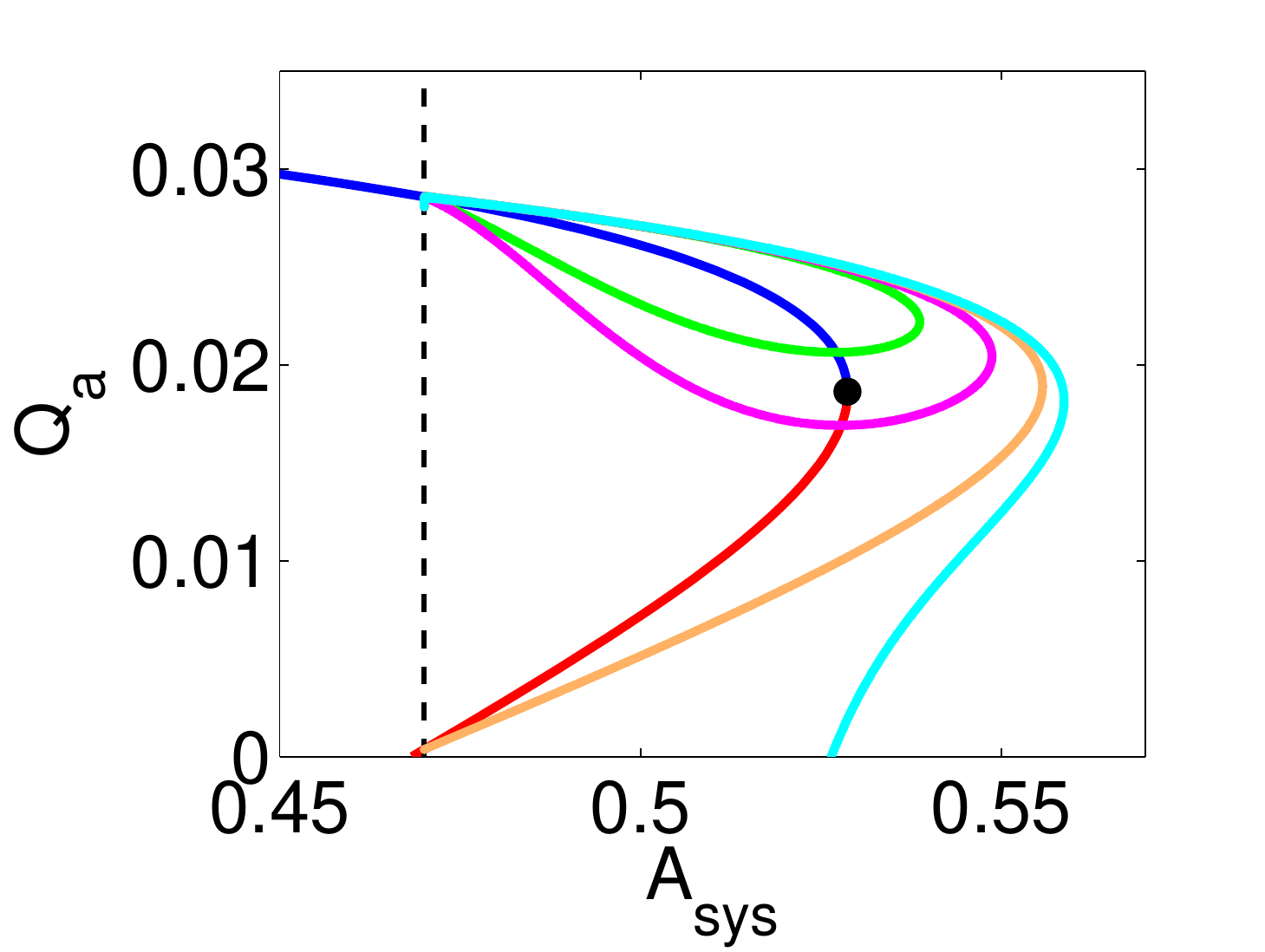}}
        ~ 
        \caption{\prrsa{(\ref{Albedo time profiles:intro}) Time
            profiles of planetary albedo forcing $\Asys(t)$ with
            maximal exceedances $R = 0.01$ (green), $R = 0.02$ (pink), $R \approx 0.027$ (light brown)
            and $R = 0.03$ (bright blue), starting from present day
            albedo $\Asys^\infty=0.47$. Horizontal dotted line
            indicates location $\Asys^b$ of the critical threshold for
            fixed albedo (fold bifurcation). (\ref{Humidity time
              profiles:intro}) Time profiles of the specific humidity
            $Q_a$.  (\ref{Monsoon tajectories:intro}) Trajectories
            from (\ref{Albedo time profiles:intro}), (\ref{Humidity
              time profiles:intro}) in the $(A_{\mathrm{sys}},Q_a)$ -
            plane and $\Asys$-dependent location of equilibria for
            fixed albedo $\Asys$: upper branch stable (blue) and lower
            branch is unstable (red). Present day $\Asys^\infty$ is
            marked by dashed line, critical equilibrium value marked
            by black dot. Underlying equations are given in \eqref{qa dot}--\eqref{Ta dot} and \eqref{Albedo forcing} in Section~\ref{sec:Monsoon}.}}\label{Monsoon deterministic tipping:intro}
\end{figure}
Figure~\ref{Monsoon deterministic tipping:intro} demonstrates the
general effect, described for general systems in
Section~\ref{sec:Theory}. The specific graphs have been computed for
Zickfeld's model for the Indian summer monsoon, a model with two
time-dependent variables depending on the planetary albedo $\Asys$ as
a parameter, given in Section~\ref{sec:Monsoon}. The current planetary
albedo is estimated as $\Asys^\infty=0.47$ \cite{zickfeld2005indian},
where the monsoon is at a stable equilibrium with specific humidity
$Q_a=0.03$.  The tipping threshold for the shutdown of the monsoon is
at $\Asys^b\approx0.53$ according to our simple model (dashed line in
Figure~\ref{Albedo time profiles:intro}). For sufficiently slow
monotone increases of the albedo $\Asys$, the value $\Asys^b$ is the
tipping point, where \js{in the model} the monsoon shuts down. However,
Figure~\ref{Albedo time profiles:intro} shows four scenarios for time
profiles of the albedo, $\Asys(t)$, where it temporarily crosses the
tipping threshold $\Asys^b$ but then returns back to its present day
value $\Asys^\infty$. Whether the monsoon shuts down depends on how
far the albedo exceeds the tipping threshold maximally (the
\emph{exceedance amplitude} $R=\max\Asys(t)-\Asys^b$), and for how
long (\emph{the exceedance time} $t_e$), for each of these
scenarios. Our general deterministic result applies for the case that
the change of the albedo is \emph{slow} compared to the internal time
scales of the monsoon model, such that we may introduce a small
parameter $\epsilon$ expressing this ratio of time scales.  The
criterion whether the monsoon avoids shutdown depends only on the
parameters $R$ and $t_e$ of the albedo time profile and one system
dependent quantity $d^b$ for small $\epsilon$. To avoid shutdown the
exceedance amplitude $R$ has to be of order $\epsilon$ and
\begin{align}\label{intro:bound}
  d^b\,R\,t_e^2&\leq 16+O(\epsilon)\mbox{.}
\end{align} Criterion \eqref{intro:bound} is general for
scenarios when one temporarily crosses a critical threshold
(mathematically a fold bifurcation), the parameter time profile can be
approximated by a parabola close to the threshold, and the ratio
$\epsilon$ of parameter drift speed to internal time scales is
sufficiently small.

The consequences of our general result are shown in Figures
\ref{Humidity time profiles:intro} and \ref{Monsoon
  tajectories:intro}. One of the time-dependent variables in the monsoon
model, the specific humidity $Q_a(t)$, experiences a temporary drop in
the first two scenarios, which meet criterion \eqref{intro:bound}
(green: $R=0.01$, $t_e=8$, and pink: $R=0.02$, $t_e=11$). In
contrast, $Q_a(t)$ drops far below its stable equilibrium value for
the fourth scenario, which does not satisfy \eqref{intro:bound} (bright
blue: $R=0.03$, $t_e=13.3$). Whilst the third scenario shows the solution for $Q_a(t)$ at the boundary of criterion \eqref{intro:bound} (light brown: $R\approx0.027$, $t_e=12.7$).   Figure~\ref{Monsoon tajectories:intro}
shows how albedo $\Asys(t)$ and specific humidity $Q_a(t)$ change in
combination. For the scenarios without tipping the trajectories form
closed curves returning to their starting (present day) equilibrium
values, while in the tipping scenario the trajectory escapes from the
region where the model is valid. Figure~\ref{Monsoon
  tajectories:intro} also contains the underlying equilibria for fixed
albedo $\Asys$. The stable equilibrium is blue, above $Q_a=0.02$,
meeting the unstable equilibrium (red) in the point that would be the
tipping point for fixed albedo (black dot at $\Asys\approx0.53$,
$Q_a\approx0.02$). The trajectory on the boundary established by equality in
\eqref{intro:bound} connects the stable and the unstable equilibrium
at the present day albedo value $0.47$ (light brown
in Figure~\ref{Monsoon tajectories:intro}).  \prrsa{We may also
  formulate a criterion equivalent to \eqref{intro:bound} which
  depends on the acceleration $\ddot{A}_{\mathrm{sys}}(t)$ of the
  planetary albedo at its maximum. This is provided by equation
  \eqref{Critical vs decay rate}.}

\paragraph{Value and estimates of  proportionality factor $d^b$}
The system dependent proportionality factor $d^b$ depends on the decay
rate \finalchange{-$\lambda(\Asys)$} toward the stable equilibrium for fixed albedo
$\Asys$ just below the critical threshold $\Asys^b$. In the ratio
\finalchange{\begin{equation}\label{intro:ddef}
  d(\Asys)=\frac{[-\lambda(\Asys)]^2}{\Asys^b-\Asys}
\end{equation}
the} decay rate \finalchange{$-\lambda(\Asys)$} is the negative of the leading
eigenvalue \finalchange{$\lambda(\Asys)$} of the linearization of the system at the stable
equilibrium, which is zero for $\Asys = \Asys^b$. Thus, both,
numerator and denominator in \eqref{intro:ddef} converge to $0$ as
$\Asys$ approaches its critical threshold $\Asys^b$ from
below. However, $d(\Asys)$ has a well defined limit \js{from the left}, defining $d^b$:
$d(\Asys)\to d^b$ as $\Asys\to\Asys^b$ \js{from below}.
}


In real-world applications there is often no direct access to an
underlying model. Instead, only time series output data, disturbed by
random fluctuations, may be available, such as proxies for the
temperature and $CO_2$ \citep{ditlevsen2010tipping} in palaeo-climate
records. \prrsa{For a constant parameter $\Asys^c<\Asys^b$ and close to the
  stable equilibrium, the system, when subjected to small white-noise
  disturbances, acts like a linear Ornstein-Uhlenbeck process where
  the decay rate \finalchange{$-\lambda$} is approximately related to the
  time-$\Delta t$ autocorrelation $a$ via \citep{aalen2008survival}
\begin{equation*}
\finalchange{-\lambda(\Asys^c)} \approx \dfrac{1-a}{\Delta t}\mbox{.} 
\end{equation*}
Thus, we} may estimate the quantity \prrsa{$d(\Asys^c)$ for some fixed
$\Asys^c<\Asys^b$ (but still $\Asys^c\approx \Asys^b$) by
\finalchange{$d(\Asys^c)\approx (1-a)^2/[(\Delta t)^2(\Asys^b-\Asys^c)]$}, where $a$
is the time $\Delta t$ autocorrelation of an output time series (e.g.,
specific humidity $Q_a(t)$ or atmospheric temperature $T_a(t)$), observed
for the fixed parameter $\Asys^c$. Then we use this as an
approximation for $d^b$ since $d(\Asys^c)$ approaches $d^b$ for
$\Asys^c\to \Asys^b$. Inserting this approximation into
\eqref{intro:bound}, we obtain a simplified dimensionless criterion
that a parameter change policy $\Asys(t)$ avoids tipping if
\finalchange{\begin{displaymath}
  \frac{(1-a)^2}{\Asys^b-\Asys^c}\left[\max_t\Asys(t)-\Asys^b\right]N_e^2\leq 16+O(\epsilon)\mbox{.}
\end{displaymath}
In} this criterion the autocorrelation was measured at $\Asys^c$
(for example, present-day value $\Asys^\infty$), and $N_e$ is the number of time units above threshold
$\Asys^b$, using the same time unit as for measuring the
autocorrelation $a$ (so $N_e=t_e/\Delta t$).} The autocorrelation
$a$ and the variance of output time series are expected to increase when
the parameter $\Asys^c$ approaches the tipping threshold $\Asys^b$
from below. This has motivated extensive studies in field data (such
as palaeo climate records or lake sediments), investigating whether
autocorrelation and variance act as \emph{early-warning indicators} of
tipping
\citep{scheffer2009early,dakos2008slowing,lenton2011early,scheffer2012anticipating}. See
also \citet{ritchie2015early} for a study on the behavior of
early-warning indicators when the parameter is changed at higher
speed, causing rate-induced tipping, \prrsa{and
  \citet{dakos2008slowing} for methods on calculating autocorrelation
  and variance}.

\paragraph*{Probabilistic result} (see Section~\ref{sec:Adding noise})
If the system is subject to small white-noise disturbances of
variance $2D$, tipping may occur with positive probability even if the
albedo time profile $\Asys(t)$ never exceeds the critical threshold
$\Asys^b$ (that is, $R<0$). Thus, the exceedance time
$t_e(\Asys^\mathrm{th})$, measuring the time the parameter $\Asys(t)$
spends above a fixed threshold value $\Asys^\mathrm{th}<\Asys^b$
becomes a relevant parameter.  We find that the level curves of
constant probability of escape follow the inverse square in parts of
the parameter space if the time scale ratio $\epsilon$ and the noise
variance have the relation $\epsilon\sim D^{2/3}$. If $\epsilon\gg
D^{2/3}$, the probabilistic result reverts to the deterministic case,
while for $\epsilon\ll D^{2/3}$ the tipping probability is close to
$1$.
We provide a numerically computed graph of escape probabilities that
is accurate in the limit $\epsilon\to0$ and $D^{2/3}\sim \epsilon$
(fitting coefficients are in the Supplementary Material). We also
provide approximations for the escape probabilities in several
limiting cases of the provided graph.

In Section \ref{sec:prob:monsoon} we will illustrate the escape
probability estimates for the monsoon model with additive
noise. The Supplementary Material provides detailed expressions
  for the projection of a general $n$-dimensional system onto a scalar ODE,
  and for the approximations of the escape probability in the presence
  of noise.

\section{Critical distance and time over threshold before tipping}
\label{sec:Theory}

\paragraph{Assumptions on the dynamical system}
We consider an $n$-dimensional system of ordinary differential equations
(ODEs) with a scalar output $y_o$
\begin{equation}
  \begin{aligned}
    \dot{\y}(t) &= f(\y(t),q(t))\mbox{,}& \y(t)&\in\mathbb{R}^n\mbox{,} \, q(t)\in\mathbb{R}\\
    y_o(t)&=\w^T\y(t)\mbox{,}& y_o(t)&\in\mathbb{R}\mbox{,} \,
    \w\in\mathbb{R}^n\mbox{.}
  \end{aligned}
\label{N ODE}
\end{equation}
that has a fold (saddle-node) bifurcation for constant $q$ at $(\y,q)
= (\y^b,q^b)$.  Specifically, we make the following assumptions
  ((S1)-(S4) define the fold bifurcation \citep{K04}):
\begin{compactenum}
\renewcommand\labelenumi{(\theenumi)}
\renewcommand{\theenumi}{S\arabic{enumi}}
\item\label{ass:singular} the linearization $A_1=\partial_1f(\y^b,q^b)$ is
singular and has a single right nullvector $\v_0$ and a single left
nullvector $\w_0$ ($A_1\v_0=0$, $\w_0^TA_1=0$), scaled such that
$\w_0^T\v_0=1$;
\item all other eigenvalues of $A_1$ have negative real part: one of the branches of the fold is stable;
\item $a_0:=\w_0^T\partial_2f(\y^b,q^b)\neq0$: changing the
  parameter $q$ crosses the fold  transversally;
\item
  $\kappa:=\frac{1}{2a_0}\w_0^T\partial_1^2f(\y^b,q^b)\v_0^2\neq0$:
  only one node and one saddle collide in the fold;
\item\label{ass:output} $\w^T\v_0\neq0$: one can observe the dynamics in the critical
  direction $\v_0$ through the output $y_o$; thus we scale $\v_0$ such
  that $\w^T\v_0=1$.
\end{compactenum}
\prrsa{We use the convention that $\partial_1^kf(\y,q)$ and
  $\partial_2^kf(\y,q)$ refer to the $k$th order partial derivatives
  of $f$ with respect to $\y$ and $q$ respectively.}  Without loss of
generality we assume that the stable equilibrium involved in the fold
exists for $q<q^b$ and the stable \prrsa{equilibrium has} output
$y_o<y_o^b:=\w^T\y^b$, such that the signs of $a_0$ and $\kappa$ are
positive:
\begin{align*}
  a_0&>0\mbox{,}&\kappa&>0\mbox{.}
\end{align*}
Otherwise, we may change the sign of the considered output projection
$\w$ or parameter $q$.  

\paragraph{Assumptions on the forcing}
\prrsa{We assume that the dependence of the parameter on time, $q(t)$,
  is slow, differentiable sufficiently often, and that $q(t)$ reaches
  a maximum, which we can assume without loss of generality to be at
  time $t=0$. Thus, $q(t)$ can be split into a constant part (equal to
  $q^{\max}$) and a time-dependent part $q_h$ with a time dependence
  of the form $\epsilon t$ and a small parameter $\epsilon$: $q(t) =
  q^{\max} + q_h(\epsilon t)$, where, by our assumptions, $q_h(0)=\dot
  q_h(0)=0$. The slowness of the parameter change implies that for
  $q^{\max}>q^b$ the system will tip for sufficiently small
  $\epsilon$, unless the overshoot of the forcing beyond the
  bifurcation value $q^b$, the maximal exceedance $q^{\max}-q^b$, is
  small. Thus, we may introduce a rescaled measure $R_0$ of the
  maximal exceedance such that $q^{\max} = q^b + \epsilon R_0$ (the
  analysis below will show that for small $\epsilon$ the boundary for
  tipping occurs for $R_0$ of order $1$). Furthermore, we assume that
  the parameter \prrsa{forcing} $q(t)$ reaches a regular maximum at
  time $t=0$, such that $R_2:=-\frac{1}{2}q_h''(0)>0$ (using
  $(\cdot)'$ to indicate the derivative of a single-argument function
  with respect to its argument). In summary, we assume that the
  parameter forcing is of the form
  \begin{equation}
    q(t) = q^b + \epsilon R_0 + q_h(\epsilon t)\mbox{,}
    \label{q expansion}
  \end{equation}
  where $\epsilon$ is small and
\begin{compactenum}\renewcommand\labelenumi{(\theenumi)}
\renewcommand{\theenumi}{P\arabic{enumi}}
\item\label{ass:qbif} $q_h(0)=0$: the parameter $q(t)$ reaches the
  value $q^{\max}=q^b+\epsilon R_0$ (without loss of generality) at time $0$;
\item\label{ass:qmax} $q_h'(0)=0$: the value $q^{\max}$ is a critical
  point of the parameter dependence, making the encounter of the fold
  at $q^b$ \emph{non-transversal} for $R_0=0$;
\item\label{ass:qcurve} $R_2:=-\frac{1}{2}q_h''(0)>0$: the parameter
  \prrsa{forcing} has a regular maximum at time $t=0$.
\end{compactenum}} Assumption \ref{ass:qmax} implies that we are
studying the vicinity of a degeneracy.  Commonly, one assumes that the
crossing of the bifurcation is transversal (that is, $q_h'(0)>0$,
making \ref{ass:qcurve} unnecessary)
\citep{kuehn2015multiple,neishtadt1987persistence,neishtadt1988persistence,baesens1991slow,baer1989slow}. The
case $q_h'(0)>0$ is called a slow (in comparison to the response time
of the system) transversal passage through a fold bifurcation.  For
the transversal case it has been shown
\citep{berglund2006noise,majumdar2013transitions} that solutions track
the stable equilibrium branch for $t<0$ at a distance of order
$\epsilon^{1/3}$ and that solutions gain distance of order $1$ from
the equilibrium branch with a delay of order $\epsilon^{2/3}$ such
that $y(t)-y^b\sim O(1)$ for times $t\sim \epsilon^{2/3}$.
\prnew{\citet{li2016time} determine the precise time at which $\y$
  reaches order $1$, corresponding to an escape to infinity after
  rescaling, calling it the ``point of no return'' for a linear
  \prrsa{forcing}.}  The textbook by Berglund and Gentz
\citep{berglund2006noise} also derives asymptotic probabilities and
escape times for tipping in the presence of white noise in the
small-$\epsilon$-small-variance regime.


The assumptions imply that we have for fixed $q<q^b$ a branch of
stable equilibria $(\y^s(q),q)$ and a branch of unstable equilibria
\prrsa{$(\y^u(q),q)$} of \eqref{N ODE} (all satisfying
$0=f(\y^{s,u}(q),q)$), which meet in the fold at parameter $q^b$ in
point $\y^b$. We expect a solution $\y(t)$, starting from close to
  $\y^s(q(t_0))$ with $t_0<0$ to follow the stable branch closely for
  sufficiently small $\epsilon$ until we reach the vicinity of $\y^b$
  at time $t<0$ of order $1$.

In the vicinity of the fold $(\y^b,q^b)$, we may zoom in
and speed up time:
\begin{equation}\label{eq:scaling}
  \begin{aligned}
    x&:=\epsilon^{-1/2}(y_\pr{o}-y_o^b)=
    \epsilon^{-1/2}\w^T(\y-\y^b)\mbox{,}\\
    t_\mathrm{new}&:=\epsilon^{1/2}t_\mathrm{old}\mbox{.}
  \end{aligned}
\end{equation}
Then \prrsa{combined with the expansion for $q(t)$, equation \eqref{q expansion},} $x$ satisfies the  scalar differential equation \citep{ritchie2017probability}
\begin{align}
\label{ODE projection}
  \dot{x} = a_0&\left(R_0-R_2t^2+\kappa x^2\right) +
  O(\epsilon^{1/2})\mbox{.}
\end{align}
If $R_0$ is sufficiently large, the trajectory $x(t)$ will grow to
large values \prrsa{for positive $t$ and small $\epsilon$ (since
  $a_0>0$).} Thus, $\y(t)$ will leave the neighborhood of the branches
of equilibria (corresponding to
tipping). 

The system quantities $a_0$ and $\kappa$ can be estimated from
observations of the output $x$ for fixed parameter $q$ (thus,
$R_2=0$): $2\kappa$ is the curvature of the equilibrium curve as
observed through $x$ in $x=0$ \finalchange{(but also in the $(y_o,q)$ plane in $(y_o^b,q^b)$)}. \prrsa{The decay rate of \eqref{ODE
    projection} at $R_2=0$ and fixed $R_0<0$ toward the stable
  equilibrium equals} $2a_0\sqrt{-R_0\kappa}$ (recall that
$\kappa>0$). Note that this is the decay rate for the sped up time
  $t_\mathrm{new}$.

\prrsa{For small $\epsilon$}, the scalar equation \eqref{ODE projection} has
solutions that are asymptotically {$x(t)\sim -|t|\sqrt{R_2/\kappa}$} for large $t$ if (see \cite{ashwin2017parameter})
\begin{equation}\label{threshold:nf}
  R_0<\frac{1}{a_0}\sqrt{\frac{R_2}{\kappa}}+o(1)\mbox{.}
\end{equation}
\prrsa{(The term $o(1)$ stands for terms that go to $0$ for $\epsilon\to0$.)} For $\kappa
a_0^2R_0^2=R_2$ and $\epsilon=0$ the orbit $x(t)=t\sqrt{R_2/\kappa}$
is the only solution existing for all time (connecting $x\sim t$ at
$-\infty$ and $+\infty$). In the original coordinates this gives a
first-order expansion for the condition relating the maximum value of
$q(t)$ \prrsa{and its acceleration at the maximum to each other. By
  assumptions \ref{ass:qbif}--\ref{ass:qcurve} the maximum of $q(t)$
  is attained at $t=0$, $\max_tq(t)= q^b + \epsilon R_0$, and the
  acceleration of $q$ at $0$ equals $\ddot q(0)=\frac{\d^2}{\d
    t^2}q_h(\epsilon t)\vert_{t=0}=-2\epsilon^2R_2$ (by
  expansion \eqref{q expansion} of $q$). Expressed using $q(0)$ and
$\ddot q(0)$, criterion \eqref{threshold:nf} for avoiding tipping
\prrsa(including the small perturbation $\epsilon$ again)} reads
\begin{equation}
  q(0)<q_\mathrm{crit}(\epsilon):=q^b+\frac{1}{a_0}\sqrt{-\frac{\ddot
      q(0)}{2\kappa}}+o(\epsilon)\mbox{.}
\label{Critical distance}
\end{equation}
The first term added to $q^b$ is of order $\epsilon$ since $\ddot
  q(0)=\epsilon^2q_h''(0)=-\prrsa{2}\epsilon^2R_2+o(\epsilon^2)$. The combination of the 
  quantities $a_0$ and $\kappa$, needed for \eqref{Critical distance},
  $1/(a_0\sqrt{\smash[b]{2\kappa}})$, may be found via the \js{left limit}
\finalchange{\begin{align}\label{d:def}
  d^b:=\lim_{q\nearrow q^b}\frac{[-\lambda(q)]^2}{q^b-q}\mbox{,}
\end{align}
where} $\lambda(q)$ is the leading eigenvalue of the linearization of
underlying system \eqref{N ODE} toward the stable equilibrium $y^s(q)$
(or, equivalently, \finalchange{the negative of} the decay rate toward $y^s(q)$) in the original
time-scale $t_{\mathrm{old}}$ and spatial scale $\y$. Then
$d^b=4a_0^2\kappa$, such that \prrsa{the acceleration} criterion \eqref{Critical distance}
for avoiding tipping becomes
\begin{align}
  \label{Critical vs decay rate}
  q(0)<q_\mathrm{crit}(\epsilon)=q^b+\sqrt{-\frac{2\ddot q(0)}{d^b}}+o(\epsilon)\mbox{,}
\end{align}
where both, $d^b$ and $\ddot{q}(0)$, are computed in the original time
and space scale.  Thus, to establish the critical permissible distance
$\max_t q(t)-q^b$ over the threshold before tipping,
we need some estimate of the attraction rate toward the stable
equilibria near the fold. This decay rate can, for example, be
estimated through the autocorrelation in the output time series $y_o(t)$
when the system is subject to fluctuations
\citep{scheffer2009early,dakos2008slowing,ditlevsen2010tipping,lenton2012early}.

Furthermore, for every \prnew{\prrsa{forcing $q(t)$} exceeding the
  bifurcation value $q^b$} ($q(0)>q^b$) we \prnew{may, as an
  alternative to $\ddot q(0)$, consider} the exceedance time $t_e$,
the time that the parameter \prrsa{forcing} $q(t)$ spends beyond the
fold bifurcation value $q^b$. In the original time scale,
\prrsa{expanding the parameter \prrsa{forcing} $q(t)$ with respect to $t$:
\finalchange{\begin{equation*}
q(t) = q(0) + \frac{1}{2}\ddot{q}(0)t^2 + O(\epsilon t)^3,
\end{equation*}
and} establishing the times $t_{\pm}$ at which the parameter \prrsa{forcing} crosses $q^b$
\finalchange{\begin{equation*}
t_{\pm} = \pm\sqrt{\dfrac{2(q^b-q(0))}{\ddot q(0)}} + O(1),
\end{equation*}
gives}} the relationship between $t_e$ and the other \prrsa{forcing} parameters
$q(0)$ and $\ddot q(0)$ as approximately
\begin{equation}
  \label{Time duration}
  \begin{aligned}
    t_e&=\prrsa{t_+-t_-=}\sqrt{\frac{8(q^b-q(0))}{\ddot q(0)}}+O(1)
    = \sqrt{\frac{4R_0}{\epsilon R_2}}+O(1) \mbox{.}
  \end{aligned}
\end{equation}
As the second expression 
makes clear, the exceedance time $t_e$ is large (of order $\epsilon^{-1/2}$), even when
the amplitude of the exceedance $\epsilon R_0+O(\epsilon^2)$ is small. 
We can then insert relation~\eqref{Time duration} into \eqref{Critical
  vs decay rate} to eliminate $\ddot q(0)$ and establish the inverse-square
law for maximal exceedance amplitude $q(0)-q^b$ and time of exceedance $t_e$
that avoids tipping (recall that $q(0)=\max_tq(t)=q^{\max}$ and $d^b$ is given by
\eqref{d:def}):
\begin{equation}
  \label{Critical time}
  d^b\left[q^{\max}-q^b\right]t_e^2\leq 16\mbox{,\quad or,\quad }
  a_0^2\kappa\left[q^{\max}-q^b\right]t_e^2\leq 4\mbox{.}
\end{equation}
In applications the parameter \prrsa{forcing} is typically not
  given in the form of an expansion such as \eqref{q expansion}, but
  rather as a function of time $q(t)$. The quantities $d^b$,
  $q^{\max}-q^b$ and $t_e$ can be computed or estimated without
  explicitly introducing $\epsilon$. Then the above inequality is a
valid criterion for the tipping threshold if $q^{\max}-q^b$ is
  small in modulus, while $t_e$ is large and the \prrsa{forcing} is
  approximately parabolic in $[-t_e,t_e]$.

\section{Indian summer monsoon model}
\label{sec:Monsoon}

We illustrate the general deterministic criterion using a \js{conceptual} model
  for one of the recognized policy-relevant tipping elements in the
Climate System, the Indian summer monsoon
\citep{lenton2008tipping,o2013addressing}. The Indian economy and
agriculture is heavily reliant on the Indian summer monsoon
\citep{bhat2006indian} as it provides the main source of water for
India \citep{liepert2015global}. In the second half of the 20th
century, summer rainfall has decreased leading to an increasing
frequency of droughts \citep{ramanathan2005atmospheric}, reducing rice
harvests \citep{auffhammer2006integrated}.  In particular, in 2002
India experienced a major drought with a seasonal rain deficit of
$21.5\%$ \citep{bhat2006indian}, seeing an increase in suicides
amongst farmers and an estimated cost of $340$ million dollars to the
Indian government for drought relief programs
\citep{liepert2015global}. \citet{meehl2008effects} connect these
observations of decreased rainfall and increased droughts to an
already present disruption of the monsoon.

We study a model for the Indian summer monsoon \citep{zickfeld2005indian}, which contains the key driving force of the monsoon, a 
moisture-advection feedback loop \citep{levermann2009basic}. In the
summer months the land is warmer than the ocean. This temperature
difference generates winds coming off the Indian Ocean onto the
land. The winds carry moisture from the ocean which is deposited over
the land in the form of precipitation. This process releases latent
heat, meaning that the temperature over land increases. A larger
temperature difference causes stronger winds carrying more moisture
and hence the positive feedback loop is formed.

We use a \js{conceptual} model proposed by \citet{zickfeld2004modeling} and make
further simplifications, though retaining the key mechanism of the
monsoon, the positive feedback loop described above. \js{The model has two time-dependent variables,}
the specific humidity $Q_a$ and the atmospheric temperature $T_a$, are
described by the following ODEs:

\begin{align}
\label{qa dot}
\dot{Q}_a &= \dfrac{E - P + A_v}{\beta I_q}, \\
\dot{T}_a &= \dfrac{\mathcal{L}(P - E) - F_{\uparrow}^{LW,TA} + F_{\downarrow}^{SL,TA}(1 - A_{\mathrm{sys}}) + A_T}{\beta I_T}
\label{Ta dot}
\end{align}

\noindent where the terms on the right-hand side are grouped as
follows:
\begin{compactitem}
\item Evaporation $E$ (mm$/$s): Proportional to the {temperature}
  difference \del{in temperature} between the land $T_a$ and the
  Indian Ocean $T_{oc}$ and {to} the {difference between} saturated
  humidity $Q_{\mathrm{sat}}$ and specific humidity $Q_a$
\begin{equation*}
E:= E(Q_a,T_a) = C_E(T_a-T_{oc})(Q_{\mathrm{sat}}-Q_a).
\end{equation*}
\item Precipitation $P$ (mm$/$s): Proportional to the specific
  humidity
\begin{equation*}
P:= P(Q_a) = C_PQ_a.
\end{equation*}

\item Moisture advection $A_v$ (mm$/$s): Winds driven by the
  temperature difference between the land and ocean bring moisture
  from the ocean over land proportional to the humidity over the ocean
  $Q_{oc}$. Winds are reversed above a given height taking moisture
  away proportional to the humidity over land $Q_a$
\begin{equation*}
A_v:= A_v(Q_a,T_a) = (T_a-T_{oc})(C_{mo}Q_{oc}-C_{ml}Q_a).
\end{equation*}
 
\item Outgoing long-wave radiation $F_{\uparrow}^{LW,TA}$ (kg$/$s$^3$): Proportional to the temperature {of the land}
\begin{equation*}
F_{\uparrow}^{LW,TA}:= F_{\uparrow}^{LW,TA}(T_a) = C_{L1}T_a + C_{L2}.
\end{equation*}

\item Incoming short-wave radiation $F_{\downarrow}^{SL,TA}$
  (kg$/$s$^3$): Fraction of incoming solar radiation not reflected,
  proportional to $1-\Asys$, where $\Asys$ is the system planetary
  albedo.

\item Heat advection $A_T$ (kg$/$s$^3$): Winds driven by the
  temperature difference between the land and ocean bring cool air at
  a prescribed low altitude proportional to the potential temperature
  $\theta_\mathrm{oc}$ above the ocean ($\theta_{oc}$ is
  fixed). Reversed winds at a prescribed high altitude $z_h$ take warm
  air away proportional to the potential temperature above the land
  $\theta_a(Q_a,T_a)$. The potential temperature at the prescribed
  height $z_h$ is given by
  $\theta_a = T_a - (\Gamma(T_a,Q_a)-\Gamma_a)z_h$ where
  $\Gamma = \Gamma_0 + \Gamma_1(T_a-T_0)(1-\Gamma_2Q_a^2)$ (with a
  reference temperature $T_0$) is the atmospheric lapse rate and
  $\Gamma_a$ is the adiabatic lapse rate
\begin{equation*}
A_T:= A_T(Q_a,T_a) = C_H(T_a-T_{oc})(\theta_{oc} - \theta_a(Q_a,T_a)). 
\end{equation*} 
\end{compactitem}

\noindent The remaining terms are all constants. The constant ${\cal
  L}$ is the latent heat, and $\beta$ converts from seconds to decades
(the unit of time $t$ is decades). Appendix \ref{app:monsoon
  vals}, Table \ref{Table of Parameters} lists all parameters and
their values and units.

\citet{zickfeld2005indian} identified two quantities that are
influenced by human activities or subject to natural variation and
affect the stability of the monsoon. \finalchange{In the model} either an increase of the
planetary albedo $A_{\mathrm{sys}}$ or a decrease in the $CO_2$
concentration from present day values can potentially lead to a
``shutdown'' of the Indian monsoon. We will focus our analysis on the
possibility of an increase in the planetary albedo. 

The planetary albedo represents the ratio of reflected to incoming
solar radiation and can be influenced by atmospheric aerosols and
land-cover conversion \citep{zickfeld2005indian}. In
  particular, the atmospheric brown cloud haze hanging over the Indian
  subcontinent has been considered responsible for the disruption
  of the monsoon with some future projections suggesting the drought
  frequency could double within a decade
  \citep{ramanathan2005atmospheric}. This cloud haze is predominately
  made up of black carbon aerosols emanating from fossil fuel
  combustion and biomass burning, which both absorb and reflect (thus,
  increasing planetary albedo) incoming radiation
  \citep{meehl2008effects}. \citet{knopf2006multi} have previously
  performed a multi-parameter uncertainty analysis for the original
  Zickfeld model, which models the planetary albedo as a function of
  the surface albedo. They concluded that if the model is
  reliable the bifurcation point for the surface albedo is
  sufficiently far from present day values such that this point cannot
  be reached in the near future. However, brown haze is typically
  poorly captured by the functional dependency assumed in the
  models \citep{lenton2008tipping}. Thus, the pollution-driven changes
  such as brown haze may place the monsoon system closer to its
    bifurcation point than concluded by \citet{knopf2006multi}. On a
  positive note this pollution is a regional problem, so does not
  require any world-wide agreement \citep{liepert2015global} and
  therefore reversion of tipping by rapid action may be
    politically more feasible \citep{o2013addressing}. \prrsa{Note
    that, despite its name the planetary albedo is a regional property
    for the Indian subcontinent in this model.}

  \citet{zickfeld2005indian} states that the present day value of the
  planetary albedo is $A_{\mathrm{sys}}^\infty =
  0.47$. \finalchange{This is comparable} with \prnew{radiative
    transfer model output data for the period 1984-1997,} which
  suggests an averaged summer value of around $0.45$ for India
  \citep{hatzianastassiou2004long}. For fixed $\Asys$ system \eqref{qa
    dot}--\eqref{Ta dot} has a fold (saddle-node) bifurcation at
  $A_{\mathrm{sys}}^b \approx 0.53$, as shown as a black dot in the
  $(\Asys,Q_a)$-plane in Figure \ref{Monsoon tajectories:intro}.  A
  slow increase of the planetary albedo linearly beyond the fold will
  cause tipping in the monsoon model, a sudden drop in the specific
  humidity would be observed. \pr{Equation \eqref{Ta dot} highlights
    how increasing the planetary albedo affects the positive feedback
    loop outlined above. As the albedo \prnew{increases}, the change
    in temperature over land decreases, meaning a smaller temperature
    difference between the land and ocean and hence weaker winds are
    formed.}  \prnew{We will use criterion~\eqref{Critical time} to
    estimate} how long the planetary albedo \prnew{may exceed} the
  fold bifurcation parameter value $\Asys^b$ without causing a tipping
  of the Indian summer monsoon \finalchange{in the model}. \js{While
    the general tipping criterion \eqref{Critical time} could in principle be used
    to guide policy, model \eqref{qa dot},\,\eqref{Ta dot} is only
    conceptual such that the precise figures for $\Asys(t)$ serve as an
    illustration of the general result, rather than as policy
    guidelines.} We assume a change of planetary albedo as a temporary
  increase from a background value $\Asys^\infty$, namely
\begin{equation}
  \Asys(t) = \Asys^{\infty} + \frac{R+\Asys^b-\Asys^\infty}{\cosh(S(t_{\mathrm{end}}-2t))^2}
\label{Albedo forcing}
\end{equation}
for a time interval $[0,t_{\mathrm{end}}]$. \prrsa{The parameters $R$
  and $S$ correspond to the amount of overshoot and the speed of
  forcing respectively. They are both small and of equal magnitude,
  such that we may introduce rescaled parameters via
  \begin{align*}
R&=r_1\epsilon\mbox{,}&
S&=s_1\epsilon\mbox{.}
\end{align*}
}

Equation \eqref{Albedo forcing} describes \prnew{an increase of the}
planetary albedo towards (and, \prnew{if $R>0$,} beyond) the
\prnew{fold bifurcation value $\Asys^b$} before it returns to its
present day \prnew{(background) value $\Asys^\infty$}. \prrsa{Figure
  \ref{Albedo time profiles:intro} in the introduction shows four
  example time profiles of planetary albedo forcing \eqref{Albedo
    forcing} for illustration.} All four fix the speed of
\prrsa{forcing $S = 0.5$}, and vary $R$, the difference between
  maximal albedo and its bifurcation value $\Asys^b$. The exceedance
$R$ of the maximal albedo beyond the fold bifurcation value $\Asys^b$
(indicated by the horizontal dashed line) and the time $t_e$ the
  albedo $\Asys(t)$ spends above $\Asys^b$ are \prnew{determined} by
  $S$ and $R$ via
  \begin{displaymath}
    t_e=\prrsa{\frac{1}{S}}\sqrt{\frac{R}{\Asys^b-\Asys^\infty}}+O(R^{3/2}/\prrsa{S})\mbox{.}
  \end{displaymath}
Specifically expressing the exceedance time in quantities of order 1 and $\epsilon$
  \begin{displaymath}
    t_e=\prrsa{\frac{1}{\sqrt{\epsilon}}\frac{1}{s_1}}\sqrt{\frac{r_1}{\Asys^b-\Asys^\infty}}+O(\epsilon^{1/2})\mbox{,}
  \end{displaymath}
  we can see the exceedance time is large (of order
  $\epsilon^{-1/2}$), consistent with equation \eqref{Time
    duration}. The error is of order $\epsilon^{1/2}$ (smaller than
  the error of order 1 in \eqref{Time duration}) due to the even
  symmetry of the forcing $\Asys(-t)=\Asys(t)$. 

\prrsa{In addition to the illustrative curves shown in
  Figure~\ref{Monsoon deterministic tipping:intro}, we compare the
asymptotic approximation~\eqref{intro:bound} to the} numerically
computed critical curve separating a ``safe'' area (monsoon retained)
from the ``unsafe'' area, where escape toward shutdown occurs, in the
two-parameter plane. We choose as \prrsa{forcing} parameters the peak
exceedance beyond the fold $R$, and exceedance time $t_e$.

The critical parameters, for
which the exact (numerically computed) connecting orbit to the saddle
occurs, are shown as a \pr{blue solid} curve in Figure \ref{Monsoon
  tip region}. As discussed in Section \ref{sec:Theory}, the critical
amount by which the planetary albedo exceeds the fold value $\Asys^b$
is approximately inversely proportional to the square of the time the
planetary albedo stays above $\Asys^b$. For example, if the planetary
albedo \prnew{increases} just above the bifurcation
($R=0.005$) then the system can spend a long
time ($\sim30$ years) above the bifurcation value without shutting
down the monsoon. However, for a higher maximum of $\Asys$ ($R = 0.02$
above the bifurcation value) the system can maintain the monsoon only if
the exceedance time $t_e$ is shorter ($\sim15$ years).

\begin{figure}[h!]
        \centering
                {\includegraphics[scale = 0.45]{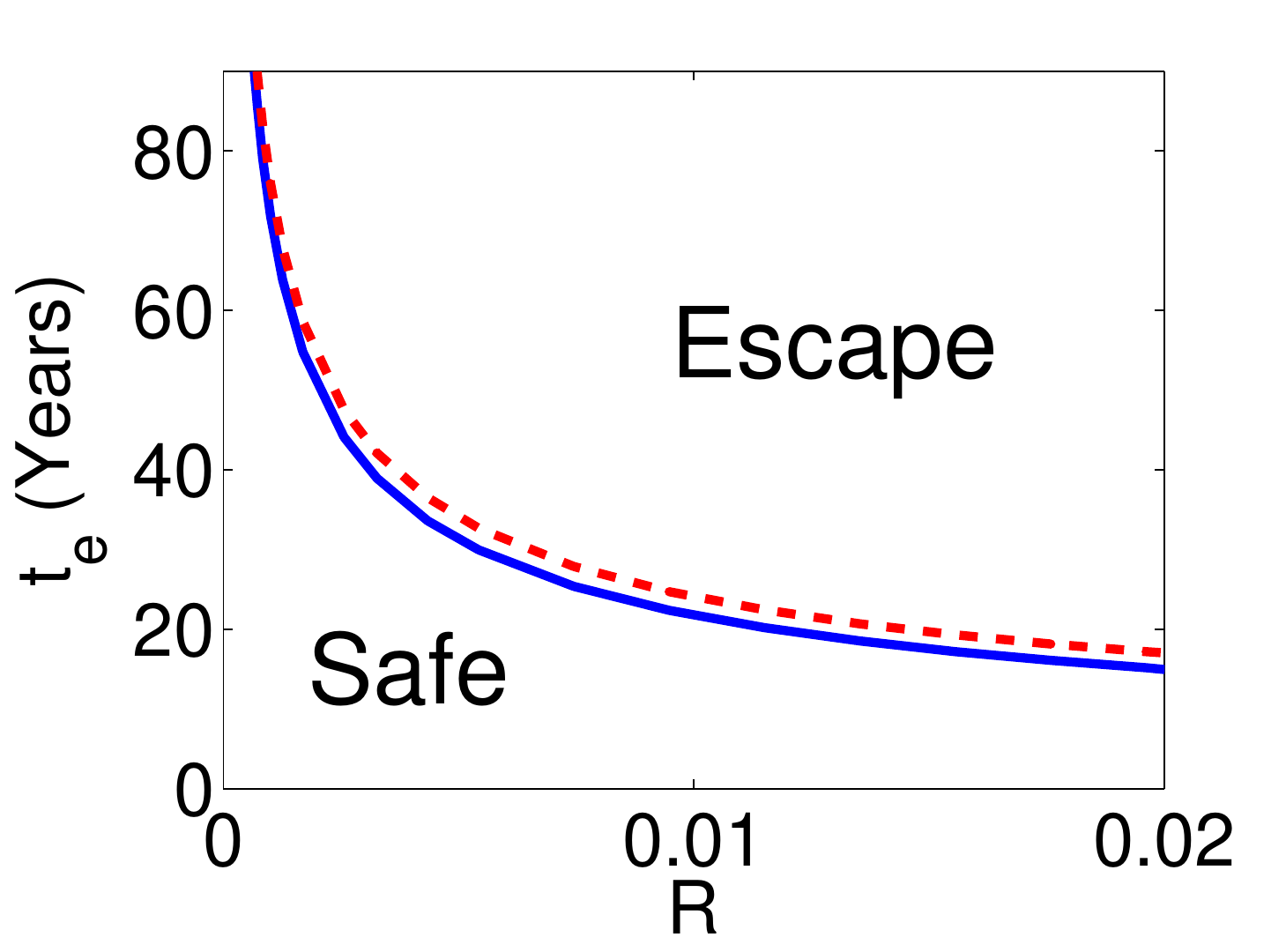}}
        ~ 
        \caption{Tipping region in the two parameter plane
          $R=\max_t\Asys(t)-\Asys^b$ (Peak distance over fold (saddle-node))
          and $t_e$ (time above fold). Safe region and escape
          region separated by the numerically calculated critical
          curve \pr{(blue solid)}. The \pr{red dashed curve} provides an
          approximation of the critical curve obtained from \pr{equation \eqref{Critical
            time}}, \changed{where $d^b = 318.36$ per decades$^2$.}}\label{Monsoon tip region}
\end{figure}

The parameter values satisfying the theoretical inequality
\eqref{intro:bound} (valid for the limit $\epsilon\to0$) are below
the red dashed curve in Figure \ref{Monsoon tip region}. The curve
gives a good approximation to the numerically calculated critical
curve. The approximation is best for small critical $R$ (peak
distance over fold) because then the system spends most time in the
region of the phase space where the second-order expansion of the
right-hand side in the fold and of the \prrsa{forcing} in its maximum are valid
(both of these were assumed in the derivation of inequality
\eqref{Critical time}, and, thus, \eqref{intro:bound}). 

\section{Probability of tipping under the influence of noise}  \label{sec:prob}
\label{sec:Adding noise}
In this section we study the probability of escape when the system
is, in addition to its parameter drift, subject to random
disturbances, which we model by adding white noise to \eqref{N
  ODE}:\prrsa{
\begin{equation}
  \label{eq:ode+noise}
  \d \y(t)=f(\y(t),q(t))\d t+\Sigma\d W_t\mbox{,}
\end{equation}
where $\Sigma$ is a $n\times\ell$ matrix of noise amplitudes, and $\d
W_t$ are the increments of $\ell$ Wiener processes}. We again consider
a \prrsa{changing parameter $q(t)$} satisfying conditions
\ref{ass:qbif}--\ref{ass:qcurve} \prrsa{and expansion \eqref{q
    expansion}}, touching a fold $(\y^b,q^b)$ of the deterministic
part satisfying conditions
\ref{ass:singular}--\ref{ass:output}. \prrsa{Similar to the analysis
  of the deterministic case in Section~\ref{sec:Theory}, the dynamics
  of \eqref{eq:ode+noise} can be studied near the fold $(\y^b,q^b)$ by
  projecting it onto its center direction using the right ($\v_0$) and
  left ($\w_0$) eigenvectors of $\partial_1f(\y^b,q^b)$.  We focus on
  the case where \finalchange{after this projection the noise amplitude is sufficiently small
  (say, of order $\sigma\ll1$)} such that escape is
  unlikely at times when the parameter \prrsa{forcing} $q(t)$ is away
  from its maximum $q^{\max}$ (for $|t|\gg1/\epsilon$). The behavior
  of the noise-disturbed system depends on the asymptotic relation
  between small noise variance of order $\sigma^2$ and small parameter
  drift speed $\epsilon$.

Two sections in the Supplementary Material derive in detail that for
$\sigma^2\ll\epsilon^{3/2}$ the probabilistic case reverts to the
deterministic case studied in Section~\ref{sec:Theory}, while for
$\sigma^2 \gg \epsilon^{3/2}$ the probability of \prrsa{tipping} will
approach $1$ for $\sigma\to0$ (and, hence, $\epsilon\to0$). Thus, for
a non-trivial limit of small noise variance $\sigma^2$ and parameter
drift speed $\epsilon$, we require a scaling of
$\sigma^2=\epsilon^{3/2}$.  The argument in the Supplementary Material
follows the textbook of \citet{berglund2006noise}, exploiting that
close to the fold the decay rate in the center direction $\v_0$ is
much smaller than the decay rates in the stable directions $\y_s$
(defined by $\w_0^T\y_s=0$), such that the coupling between stable and
center directions is small. One underlying assumption is that the
scaled projection $\sigma^{-1}\w_0^T\Sigma$ of the noise onto the
center direction is of order $1$ (that is, it is not much smaller than
the scaled projection $\sigma^{-1}[I-v_0\w_0^T]\Sigma$ onto the stable
directions).}

\prrsa{Consequently, the projection of $\Sigma$ by $\w$ onto the
  scalar output, after the rescaling \eqref{eq:scaling} to the
  zoomed-in output $x=\epsilon^{1/2}\w^T(\y-\y^b)$ and sped up time
  $\epsilon^{-1/2}t$ (see \eqref{eq:scaling}), has the variance
($\Delta=\Sigma\Sigma^T$)}
\begin{align}
  \label{eq:noiseconversion}
  2D:=2\epsilon^{-3/2}\w_0^T\Delta\w_0\mbox{,}
\end{align}
which is of order $1$.

\prrsa{Furthermore, if the \finalchange{the matrix of noise amplitudes} $\Sigma$
  in \eqref{eq:ode+noise} depends on the state $\y$, then after the
  rescaling \eqref{eq:scaling} the dependence of $D$ on $x$ becomes
  weak, of order $\epsilon^{1/2}$, such that we neglect it to leading order.}

In the limit for $\epsilon\to0$, the projected equation
\eqref{eq:ode+noise} becomes the scalar stochastic differential
equation (SDE)
\begin{align}\label{sde:unscaled}
  \d x&=a_0[R_0-R_2t^2+\kappa x^2]\d t+\sqrt{2D}\d W_t\mbox{}
\end{align}
(recall that $a_0>0$ and $\kappa>0$ without loss of generality),
\prrsa{starting from $x(t_0)\leq0$ and $t_0<0$ (see Supplementary
  Material).} By further rescaling $x$ and time and introducing
correspondingly rescaled versions of the parameters $R_0$ and $R_2$,
\pr{\begin{equation}\label{nf:scaling}
  \begin{aligned}
    x_\mathrm{new}&=\frac{(a_0\kappa)^{1/3}}{D^{1/3}}x_\mathrm{old}\mbox{,} &
    t_\mathrm{new}&=D^{1/3}(a_0\kappa)^{2/3}t_\mathrm{old}\mbox{,}&
    p_0&=\frac{a_0^{2/3}R_0}{D^{2/3}\kappa^{1/3}}\mbox{,} &
    p_2&=\frac{R_2}{D^{4/3}\kappa^{5/3}a_0^{2/3}}\mbox{,}
  \end{aligned}
\end{equation}}
we may simplify \eqref{sde:unscaled} to a SDE
\begin{align}\label{sde:scaled}
  \d x&=[p_0-p_2t^2+x^2]\d t+\sqrt{2}\d W_t\mbox{}
\end{align}
with unit noise amplitude and nonlinear coefficient, and the two
parameters $p_0\in\R$ and $p_2>0$.
The lines $x=\sqrt{p_2}t$ for $t\ll-1$ and $x=-\sqrt{p_2}t$ for
$t\gg1$ are stable slow manifolds of the deterministic part of
\eqref{sde:scaled}. Thus, the density of $x$ at some fixed large time
$t=-T_0$ is nearly independent from the initial density for $t\ll-T_0$
(conditional on no escape occurring before $t=-T_0$). Thus, we can
compute numerically the probability of escape by solving the
Fokker-Planck equation (FPE) \citep{risken2012fokker} for the density $u(x,t)$ of $x$
\begin{equation}
\partial_t u(x,t)= \partial_x^2 u(x,t)-\partial_x[(p_0-p_2t^2+x^2)u(x,t)]
\label{fullFPE}
\end{equation}
with Dirichlet boundary conditions
$u(x_\mathrm{bd},t)=u(-x_\mathrm{bd},t)=0$ from $t=-T_0$ to $t=+T_0$,
starting from an arbitrary density concentrated in the region
$\{x\leq0\}$ and a sufficiently large $T_0$. The resulting escape
probability $P_\mathrm{esc}(p_0,p_2)$ is then (approximately for large
$T_0$ and large $x_\mathrm{bd}$) given by
\begin{equation}
  \label{eq:escp}
  P_\mathrm{esc}(p_0,p_2)=1-\int_{-x_\mathrm{bd}}^{x_\mathrm{bd}}u(x,T_0)\d x\mbox{.}
\end{equation}
The result (using \texttt{chebfun} \citep{Driscoll2014}) is shown in
Figure~\ref{fig:prob_R_te}. \prrsa{Since in the probabilistic
  scenario, tipping is possible also for \prrsa{forcing}s $p(t) = p_0
  - p_2t^2$ that do not exceed the critical value $p_0^b = 0$, it is
  useful to also consider other (especially lower) thresholds than
  $p_0^b$ \prrsa{(indicated by the white dashed line)} for measuring
  exceedance amplitudes and times. In Figure \ref{fig:prob_R_te} we
  choose $p_0^\mathrm{th}=-1$, and, hence, consider the parameter
  plane $(R^{(-1)},t_e^{(-1)})$ where
  \begin{align*}
    R^{(-1)} &= p_0 -
 p_0^\mathrm{th}\in[0,3]\mbox{,\quad and}&   
 t_e^{(-1)} &= 2\sqrt{(p_0-p_0^\mathrm{th})/p_2}\in[1,5]
  \end{align*}
  are the maximal distance over the threshold $p_0^\mathrm{th}$ and the
  exceedance time spent above $p_0^\mathrm{th}$.
  
  \begin{figure}[h!]
  \centering
  \includegraphics[scale=0.5]{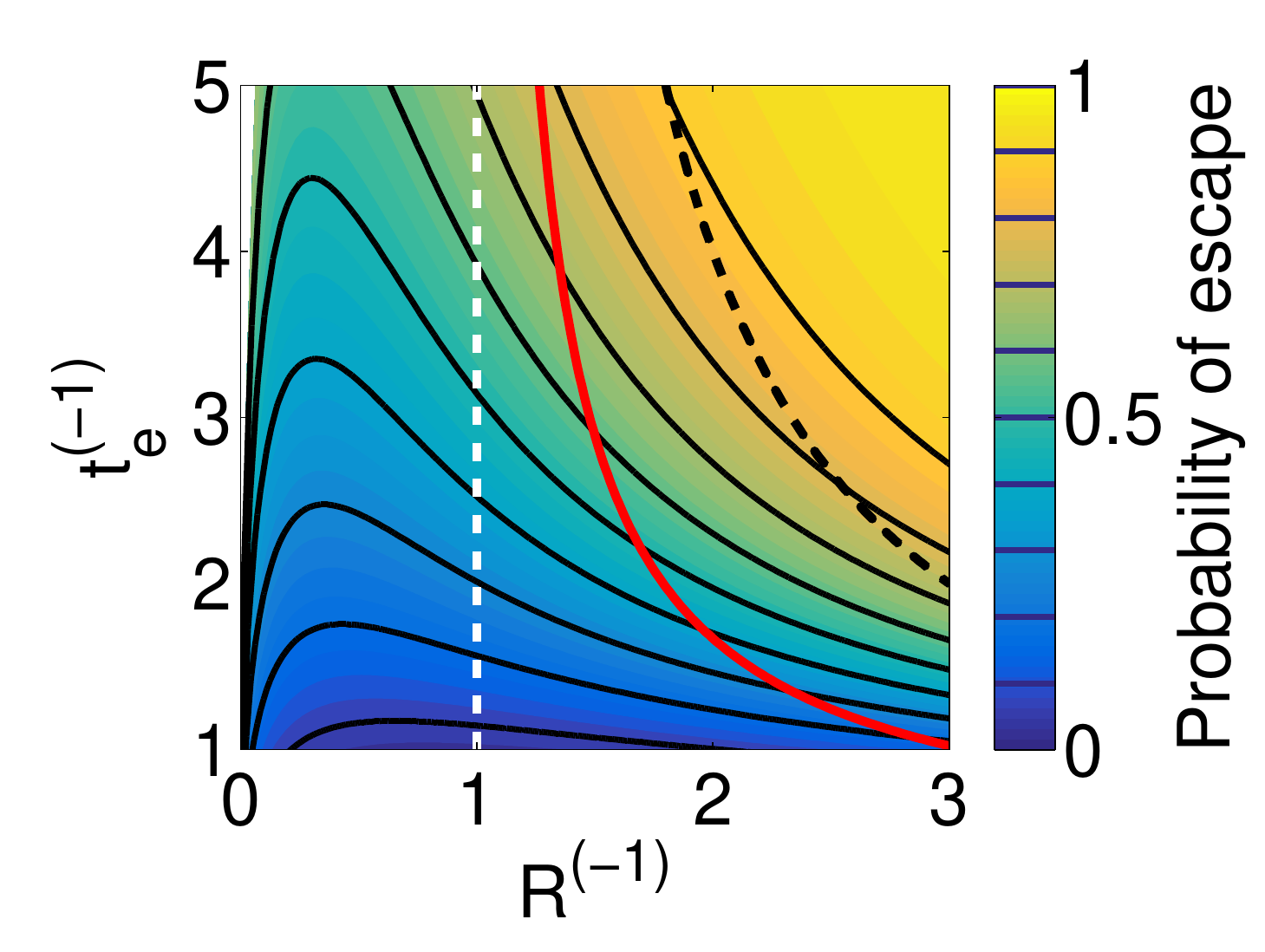}
  \caption{Probability of escape, $P_\mathrm{esc}$ of $x\to+\infty$ in
    \eqref{sde:scaled} for parameter values \prrsa{$R^{(-1)} = p_0 + 1$ and $t_e^{(-1)} = 2\sqrt{(p_0+1)/p_2}$ (exceedance amplitudes and times beyond the arbitrary threshold $p_0^\mathrm{th} = -1$)}. Parameters for
    FPE~\eqref{fullFPE}: domain $[-8,8]$, time interval $[-T_0,T_0]$,
    initial density $N(x_0,1)$ where $T_0=\sqrt{(x_0^2+p_0)/p_2}$ and
    $x_0=-4$. \prrsa{White dashed line indicates the static fold bifurcation, black dashed curve shows the critical parameter values for deterministic tipping and the red solid curve indicates the boundary for validity of the mode approximation.}}
  \label{fig:prob_R_te}
\end{figure}

  In the lower right corner of Figure~\ref{fig:prob_R_te}, the level
  curves of equal probability align with the inverse-square law for
  deterministic tipping (black dashed curve in Figure
  \ref{fig:prob_R_te} provides boundary for deterministic
  tipping). For fixed expansion parameters $R_0$ and $R_2$ of the
  parameter drift $q(t)$ in this corner corresponds to small noise
  variance $D$. Figure~\ref{fig:q1_q2} in
  Appendix~\ref{app:level_prob} shows a transformation of the
  $(p_0,p_2)$-plane in which varying the variance $D$ corresponds to
  moving along a straight line. For finite positive noise Figure
  \ref{fig:prob_R_te} shows the deviation from the inverse-square
  relationship. For $R^{(-1)}<1$ (meaning that the fold is not
  reached) and small $t^{(-1)}_e$ (fast shifts) the probability of
  escape is small. The probability of escape increases if the
  \prrsa{maximal forcing} exceeds the fold ($R^{(-1)}>1$) or the
  exceedance time $t^{(-1)}_e$ over the threshold increases.  The
  $(R^{(-1)},t_e^{(-1)})$ coordinates are singular at their origin
  $(0,0)$ such that all equal-probability level curves pass through
  the origin.}

\paragraph{\prrsa{Mode approximation in moving coordinates}}
While the numerical result is sufficient for some practical estimates,
the Supplementary Material provides two approximation formulas for
some regions of the \prrsa{$(R^{(-1)},t_e^{(-1)})$ plane that are
  generalizable to non-parabolic parameter changes. These
  approximations show how escape rates become exponentially small
  (similar to Kramers' escape rate approximation
  \cite{hanggi1990reaction}) such that they} provide uniform accuracy
in the region of small \prrsa{$p_2$ or large $t_e^{(-1)}$}. In
particular, the \emph{mode approximation}, first tested by
\citet{ritchie2017probability}, is valid to the left of the red curve
in Figure~\ref{fig:prob_R_te}.

The mode approximation is based on the escape rate at time $t$
provided by the leading eigenvalue $-\gamma_1(\bar x)$ of the
Fokker-Planck operator $[F_{\bar
    x}u](z):=\partial_z^2u(z)- \partial_z[(z^2+2\bar x z)u(z)]$ with
  Dirichlet boundary conditions for large $z$ at both ends. The value
  of $\bar x$ at time $t$ for parameters $p_0$ and $p_2$, $\bar
  x(t;p_0,p_2)$, is the unique solution of the deterministic part
  $\dot x=p_0-p_2t^2+x^2$ with $\bar x(t;p_0,p_2)+\sqrt{p_2}|t|\to0$
  for $|t|\to\infty$ (see the Supplementary Material for
details). This permits us to estimate escape probabilities after a
one-off fit of the leading eigenvalue $-\gamma_1(\bar x)$ as a
function of $\bar x$:
\begin{equation}
  \label{eq:mode:p}
  P_\mathrm{esc}^\mathrm{1d}\approx 1-\exp\left(\int_{-\infty}^\infty \gamma_1(\bar x(t;p_0,p_2)\d t\right)\approx 1-\exp\left(\int_{-\infty}^\infty \exp(-c_0-c_2\bar x^2(t;p_0,p_2)\d t\right)\mbox{,}
\end{equation}
where $c_0=1.01$ and $c_2=1.41$ provide a good a-priori fit of
$-\log\gamma_1(\bar x)=c_0+c_1\bar x^2$. A $4$th-order fit for
$-\log\gamma_1(\bar x)$ is given in the Supplementary Material. A
necessary condition for the accuracy of the approximation \prrsa{(see
  Supplementary Material for an explanation)} is that $\bar
x(t;p_0,p_2)<0$ for all $t$, which is the case \prrsa{left of} the red
line in Figure~\ref{fig:prob_R_te}.

In the parameter region \prrsa{to the left of} the red line in Figure
\ref{fig:prob_R_te}, \prrsa{the mode approximation has a global
  absolute error less than $0.05$ when applied to the monsoon example
  (Figure \ref{err mode} in Appendix \ref{sec:modelcomp} is a contour
  graph of the error).} 

\section{Probability of tipping in the monsoon model}
\label{sec:prob:monsoon}
We now estimate the probability of shutdown of the monsoon in the model by
projecting the system (\ref{qa dot},\ref{Ta dot}) onto a
one-dimensional output ($\w^T=\w_0^T=(-3.50,-0.99)$)
\comment{$\Leftarrow$ this correct?} and expanding it near the fold to
quadratic order (in $x$), and \prrsa{using the mode
  approximation estimate \eqref{eq:mode:p} for the escape
  probability}. If time is measured in decades, the quadratic
expansion of the monsoon model near its fold has the form $\dot x=p_f
(\Asys(t)-\Asys^b)+x_fx^2$ where $x$ is a dimensionless projection of
the state. We add white noise of variance
$2\Delta=\operatorname{diag}(0.02,6)$ \comment{$\Leftarrow$ insert!}
to the monsoon model \eqref{qa dot}--\eqref{Ta dot}, such that
\begin{equation}
  \label{eq:monsoon:nf}
  \d x=[p_f (\Asys(t)-\Asys^b)+x_fx^2]\d t+\sqrt{2D}\d W_t\mbox{}
\end{equation}
with \prrsa{the scaling factors} $p_f=115.30$, $x_f=0.69$, and
\prrsa{noise variance} $D=3.04$.  The \prrsa{albedo forcing functions}
$\Asys(t)$ are chosen identical to equation~\eqref{Albedo forcing}
\begin{equation}
  \Asys(t) = 
  \Asys^{\infty} +  \frac{R^{(0.5)}+\Asys^\mathrm{th}-\Asys^\infty}{\cosh(\prrsa{S}(t_{\mathrm{end}} - 2t))^2}\mbox{,}
\label{Albedo forcing repeat}
\end{equation}  
where $\Asys^\infty=0.47$ is the present-day value of the albedo and $[0,t_{\mathrm{end}}]$ is the time interval.      
\prrsa{For the same reasons as in the previous section,} we choose \prrsa{a threshold lower than $\Asys^b \approx 0.53$, specifically,}
$\Asys^\mathrm{th}=0.5$, \prrsa{for measuring exceedance amplitudes and times.} Hence, we consider the parameter plane
$(R^{(0.5)},t^{(0.5)}_e)$, where $R^{(0.5)}=R+\Asys^b-0.5$ is the
distance of the maximal albedo along the \prrsa{forcing} from a chosen albedo
threshold $\Asys^\mathrm{th}=0.5$, and $t^{(0.5)}_e$ is the corresponding exceedance
time above this threshold.

\begin{figure}[h!]
        \centering
       \subcaptionbox{\label{Eig prob re}}[0.32\linewidth]
                {\includegraphics[scale = 0.3]{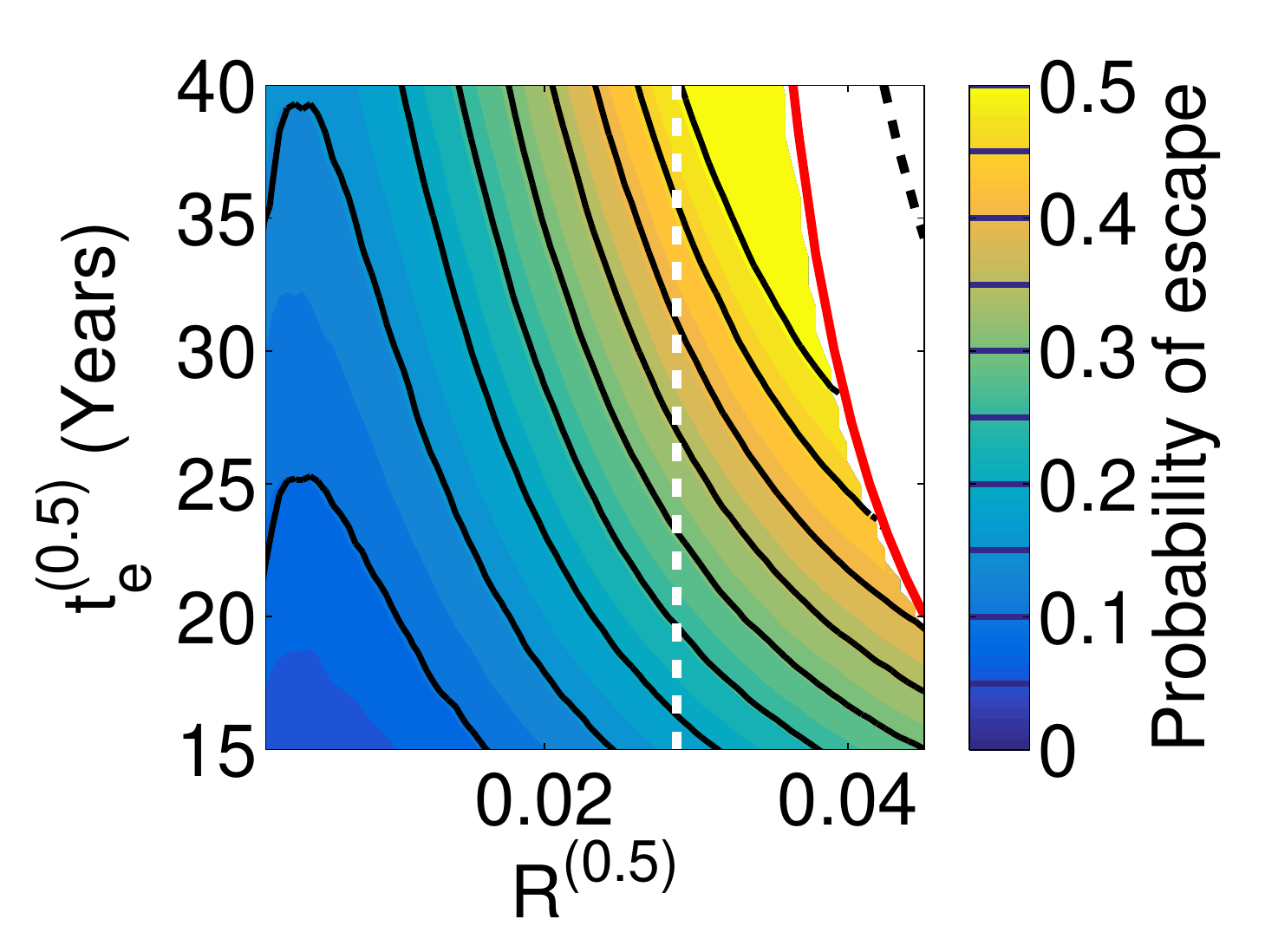}}
                \hfill
       \subcaptionbox{\label{Forcing r fixed}}[0.32\linewidth]
                {\includegraphics[scale = 0.3]{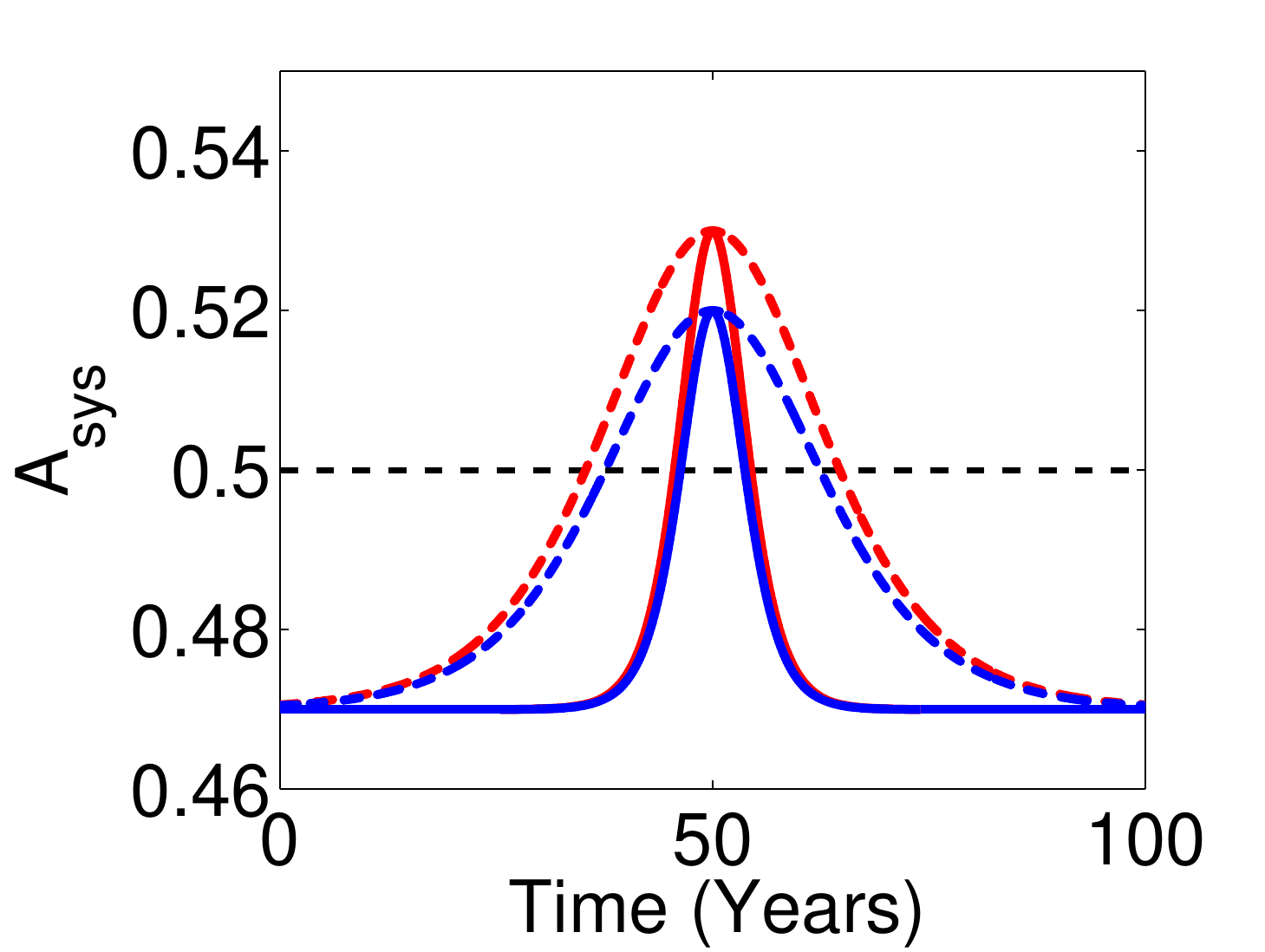}}
                \hfill
       \subcaptionbox{\label{Prob r fixed}}[0.32\linewidth]
                {\includegraphics[scale = 0.3]{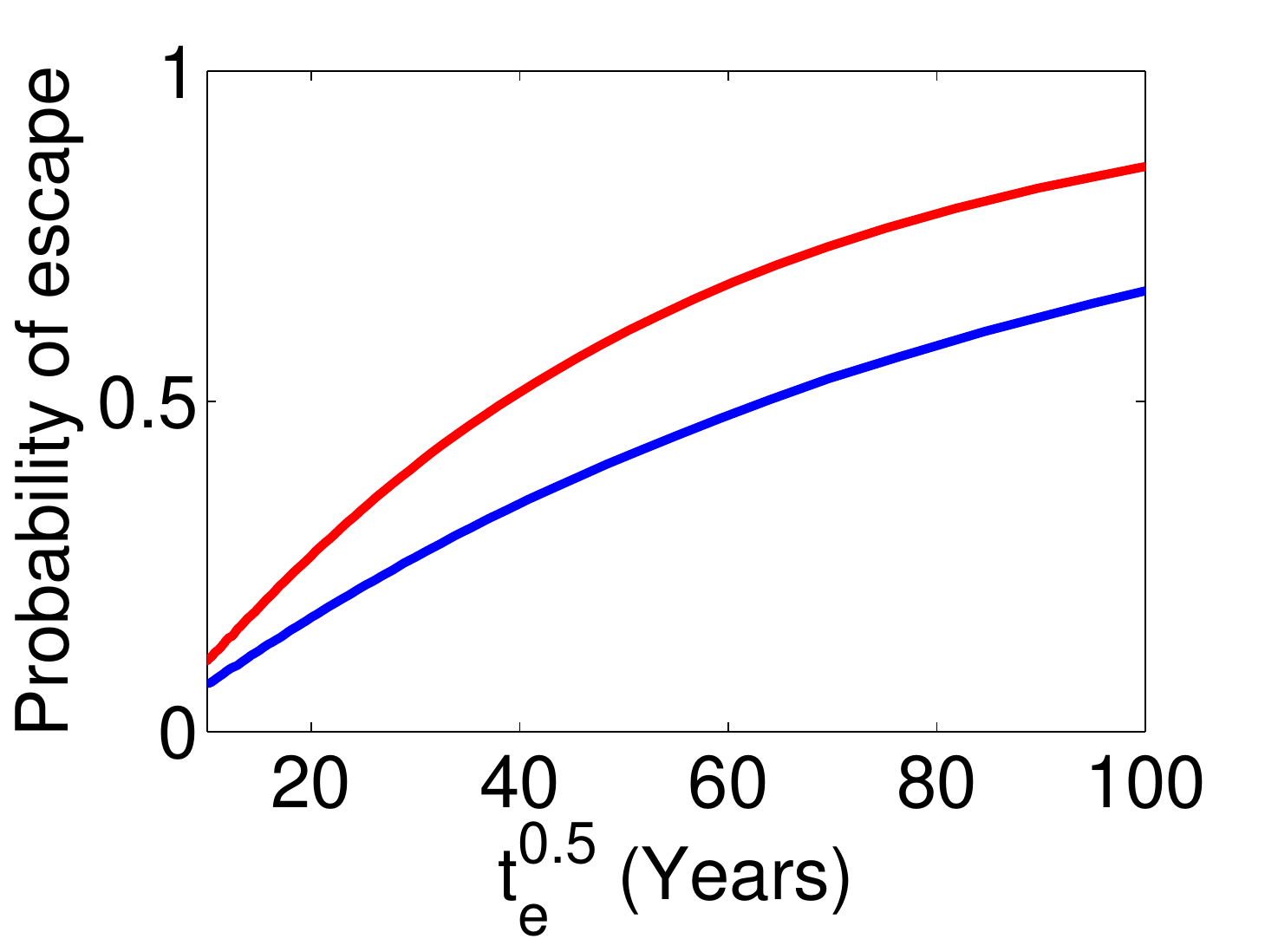}}
        ~ 
        \caption{Probability of escape: Contours in
          \pr{$(R^{(0.5)},t^{(0.5)}_e)$} - plane and cross sections
          for two different values of \pr{$R^{(0.5)}$:
            $R^{(0.5)}=0.02$ (blue) and $R^{(0.5)}=0.03$
            (red). (\ref{Eig prob re}): Contour lines are spaced at
            \pr{$0.05$} intervals. White dashed line indicates the
            static fold bifurcation, black dashed curve shows the
            critical parameter values for deterministic tipping and
            the \pr{red solid} curve indicates the boundary for
            validity of the mode approximation. (\ref{Forcing r
              fixed}) Time profiles of planetary albedo
            \prnew{scenarios} of a short \prrsa{($t^{(0.5)}_e = 7.5$,
              and $9$ years for blue solid and red solid
              respectively)} and longer \prrsa{($t^{(0.5)}_e \approx
              25$, and $30$ years for blue dashed and red dashed
              respectively)} exceedance time $t^{(0.5)}_e$ for
            \pr{each} fixed \pr{$R^{(0.5)}$} value. Horizontal black
            dotted line represents the chosen threshold value
            $\Asys^{\mathrm{th}}=0.5$ used for the definition of $R^{(0.5)}$ and
            $t_e^{(0.5)}$. (\ref{Prob r fixed}) Probability of escape
            over a range of exceedance times
            $t_e^{(0.5)}$. Parameters: {$\w=(-3.50,-0.99)^T$} such
            that $p_f=115.30$ and $x_f=0.69$ in \eqref{eq:monsoon:nf},
            $D=3.04$.}}
\end{figure}

Figure \ref{Eig prob re} shows the mode approximation
\eqref{eq:mode:p}, using the solution $\bar x$ of
$\Dot{\bar x}=p_f
(\Asys(t)-\Asys^b)+x_f\bar x^2$ with \eqref{Albedo forcing repeat}, for the
probability of escape $P_\mathrm{esc}^\mathrm{1d}$, on a
grid of points in the $(R^{(0.5)},t^{(0.5)}_e)$ - plane.  Since the
mode approximation \eqref{eq:mode:p} is valid only in the region
\prrsa{to the left of} the red line, Figure \ref{Eig prob re} leaves a
part of the $(R^{(0.5)},t_e^{(0.5)})$-plane white. \prrsa{As before,}
the vertical white dashed line positioned indicates the location of
the deterministic fold bifurcation and the black dashed curve provides
the boundary for deterministic $(D = 0)$ tipping.

\prrsa{Consistent with the previous section and as expected we observe
  that for small $R^{(0.5)}$ and small $t^{(0.5)}_e$ the probability
  of escape is small. While, increasing $R^{(0.5)}$ (to and beyond
  the fold) and/or increasing the exceedance time $t^{(0.5)}_e$ over
  the arbitrary threshold the probability of escape increases.}

Figure~\ref{Prob r fixed} shows the probability of escape for a range
of times \pr{$t^{(0.5)}_e$} {for} \pr{which $\Asys(t)$ is above the threshold $0.5$} for two
fixed maxima $\max_t\Asys(t)$. Figure~\ref{Forcing r fixed} shows
example \prnew{scenarios of $\Asys(t)$ with} two different exceedance times
for each fixed maximum.

\section{Summary} 
\label{sec:Summary}

We have investigated the scenario of forcing a system over a tipping
threshold (a fold of equilibria) for a short time. We provide simple
criteria determining whether the forced deterministic system escapes
from the family of equilibria or not. The two primary parameters of
the parameter forcing $q(t)$ are the maximum exceedance amplitude
$\max_t q(t)$ beyond the fold bifurcation parameter value $q^b$, and the
time $t_e$ for which the parameter \prrsa{forcing} has exceeded the
fold bifurcation value. The critical curve, which separates a region
of tipping and the safe region in this two-parameter plane follows an
inverse-square law: $t_e^2(\max_t q(t)-q^b)=16/d^b$. The constant $d^b$
can be determined from equilibria at parameters $q$ near the critical
value $q^b$ as the ratio of the square of the decay rate to $q-q^b$.

We used a simplified version of the Indian summer monsoon model
developed by \citet{zickfeld2004modeling} to demonstrate which
\prnew{scenarios} for changing planetary albedo in the model result in
(or avoid) tipping, matching the general theoretical predictions
(which are only accurate if one is sufficiently close to a fold)
precisely: the trade-off between exceedance amplitude and time beyond
the critical value of the albedo follows the inverse-square law with
the predicted factor $d^b$.

We also quantify the effect of random disturbances, modeled by white
additive noise, determining a probability of escape. For each chosen
threshold $q_0$ near the fold we obtain level curves of equal
probability in the parameter plane $(\max_tq(t)-q_0,t_e^0)$ of
exceedance amplitude and exceedance time for $q_0$.  These level
curves follow the inverse-square law in part, deviating from it
(expectedly) in the large-noise/slow-drift limit (for which we provide
approximation formulas in the Supplementary Material) and at the origin of the
$(\max_tq(t)-q_0,t_e^0)$.\vskip6pt

\enlargethispage{20pt}





\ethics{This work did not involve any active collection of human data.}

\dataccess{The research materials supporting this publication have been uploaded as part of the supplementary material.}

\aucontribute{P.R. derived the general mathematical results (with contributions from J.S.), processed and simplified the monsoon model, performed the numerical computations in sections 3 and 5, and led the writing of the paper. J.S. contributed to the mathematical results and performed the numerical computations in Section 4 and the fitting coefficients for approximations in the probabilistic results. O.K. carried out the comparison between approximations and simulations of the full model. All authors have contributed to the writing of the manuscript and gave final approval for publication.}

\competing{We have no competing interests}

\funding{P.R.'s research was supported by funding from the
EPSRC Grant No. EP/M008495/1, P.R. has also received funding from the NERC grant No. NE/P007880/1. J.S. gratefully acknowledges the financial support of the EPSRC via Grants No.
EP/N023544/1 and No. EP/N014391/1. J.S. has also received
funding from the European Union's Horizon 2020 research and
innovation programme under Grant Agreement No. 643073.}

\ack{We would like to thank Tim Lenton and Peter Cox for advice on the current consensus on the prospects of tipping in the Indian summer monsoon.}

%

\begin{appendix}
\setcounter{paragraph}{0}

\section{Escape probability levels in transformed
  coordinates}
\label{app:level_prob}

\begin{figure}[h!]
  \centering
  \includegraphics[scale=0.45]{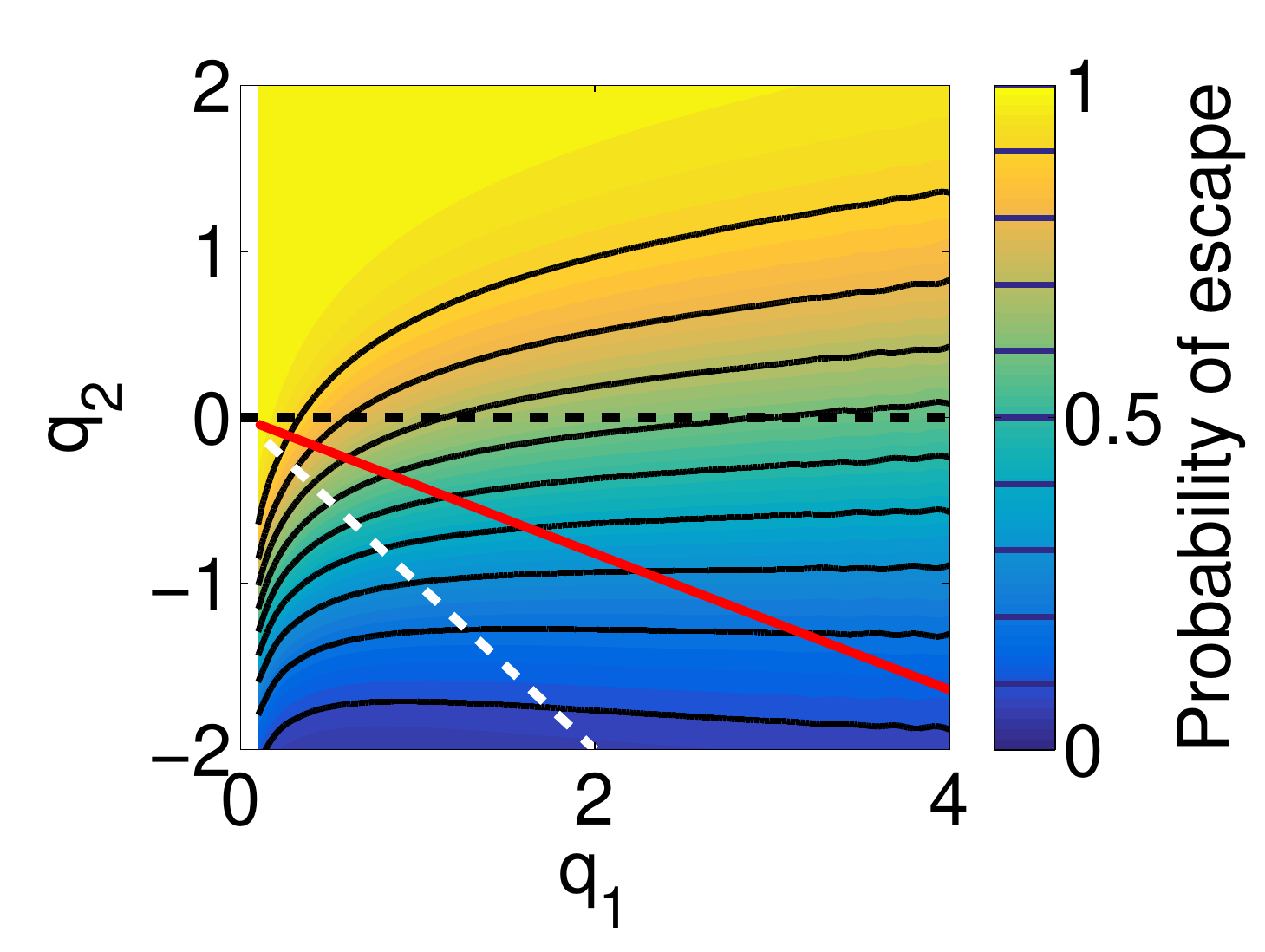}
  \caption{Transformation of the $(p_0,p_2)$-parameter plane of
    Figure~\ref{fig:prob_R_te}. for escape probability
    $P_\mathrm{esc}$ for $x\to+\infty$ in \eqref{sde:scaled} with the
    transformation $q_1=\sqrt{p_2}$ and
    $q_2=p_0-\sqrt{p_2}$. Parameters for FPE~\eqref{fullFPE}: domain
    $[-8,8]$, time interval $[-T_0,T_0]$, initial density $N(x_0,1)$
    where $T_0=\sqrt{(x_0^2+p_0)/p_2}$ and $x_0=-4$. \prrsa{White
      dashed line indicates the static fold bifurcation, black dashed
      curve shows the critical parameter values for deterministic
      tipping and the red solid line indicates the boundary for
      validity of the single-mode approximation. Computed using
      \texttt{chebfun} \citep{Driscoll2014}.}}
  \label{fig:q1_q2}
\end{figure}

\prrsa{The asymptotic behavior of the escape (tipping) probability
  for SDE \eqref{sde:scaled} becomes visible after parameter
  transformation
  $(q_1,q_2)=(\sqrt{p_2},p_0-\sqrt{p_2})\in[0.1,4]\times[-2,2]$, shown
  in Figure~\ref{fig:q1_q2}.}

When one varies the noise variance $D$ in the original
system~\eqref{sde:unscaled}, keeping the original \prrsa{forcing} parameters
$R_0$ and $R_2$ fixed, one moves along a straight line through the
origin in Figure~\ref{fig:q1_q2}. An example is shown by a
red line in Figure~\ref{fig:q1_q2}. The small-noise
limit is at the large-$q_1$ end, and the large noise (or
slow drift) limit is at the origin.

In the coordinates $(q_1,q_2)=(\sqrt{p_2},p_0-\sqrt{p_2})$ the slope
of all $P_\mathrm{esc}$ level curves of equal probability will
approach $0$ for large $q_1$ (or $\sqrt{p_2}$) such that the level
curve for $P_\mathrm{esc}=0.5$ asymptotes to the horizontal
$q_2=p_0-\sqrt{p_2}=0$ as at these parameter values the deterministic
system has its tipping threshold (see \eqref{threshold:nf} in
Section~\ref{sec:Theory}), and the limit of large $q_1$ (or $p_2$)
corresponds to the rapid drift (or small noise) limit. \prrsa{The
  $P_\mathrm{esc}$ level curves are graphs of functions $q_2(q_1)$ in
  the coordinates of Figure~\ref{fig:q1_q2}. The slopes of these
  functions, $q_2'(q_1)$ go to zero for $q_1\to\infty$. This follows
  from the fact that in the small noise limit all probabilities with
  positive distance from $0$ or $1$ must get close to the
  deterministic tipping boundary. This deterministic tipping boundary
  corresponds to the ray $\{q_2=0, q_1>0\}$. Since the sets of fixed
  forcing parameters $R_0$ and $R_2$ and varying noise variances $2D$
  correspond to rays, probabilities must converge to either $0$ or $1$
  along each ray with $q_2\neq0$. Hence, each probability contour
  level must cross all rays with non-zero slope.}

\section{Projection and approximation error for the monsoon
  model}
\label{sec:modelcomp}
\begin{figure}[ht]
  \centering
        \subcaptionbox{\label{True}}[0.32\linewidth]
                {\includegraphics[width=0.33\linewidth]{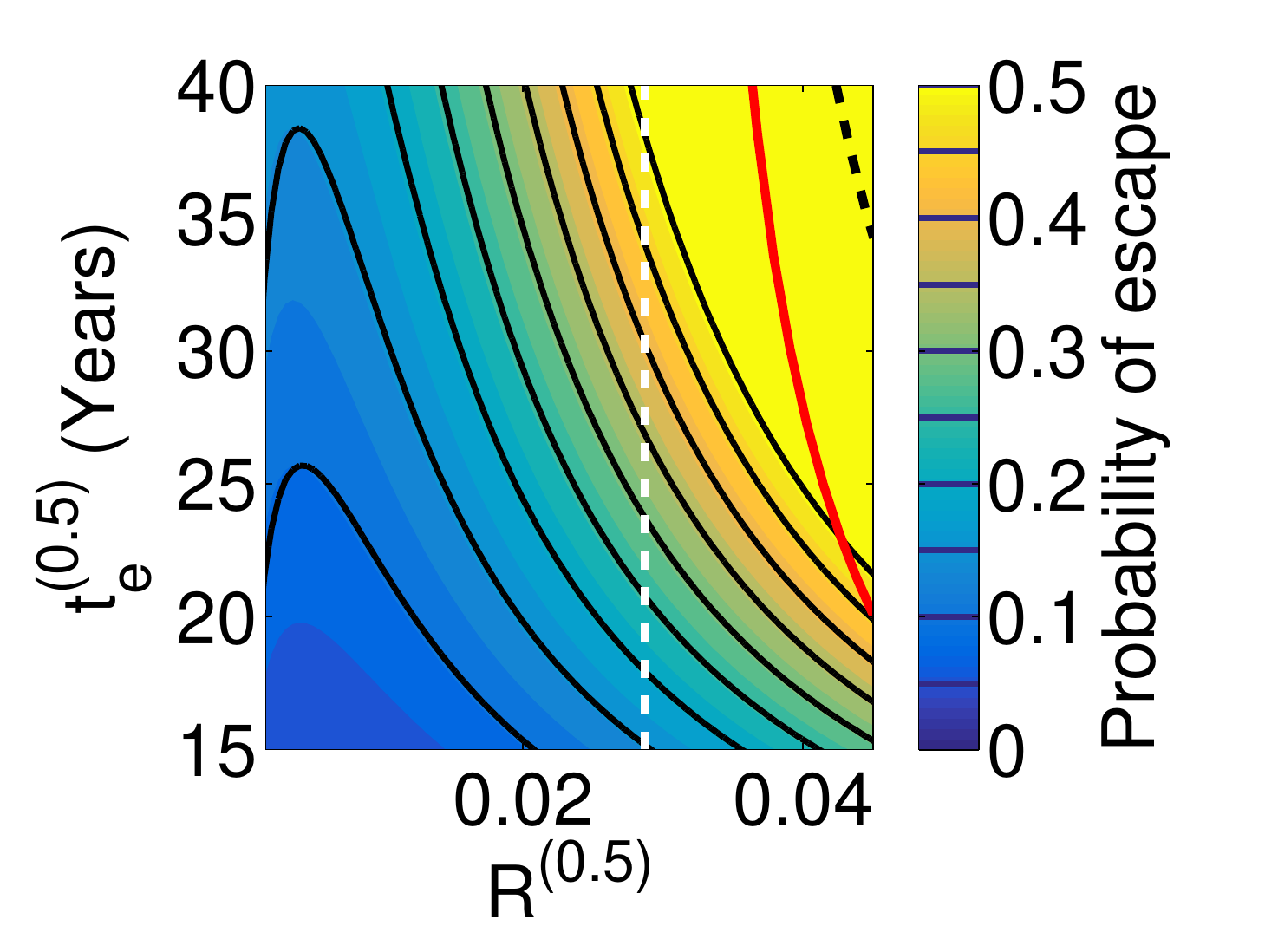}}
        \hfill
        \subcaptionbox{\label{err mode}}[0.32\linewidth]
                {\includegraphics[width=0.33\linewidth]{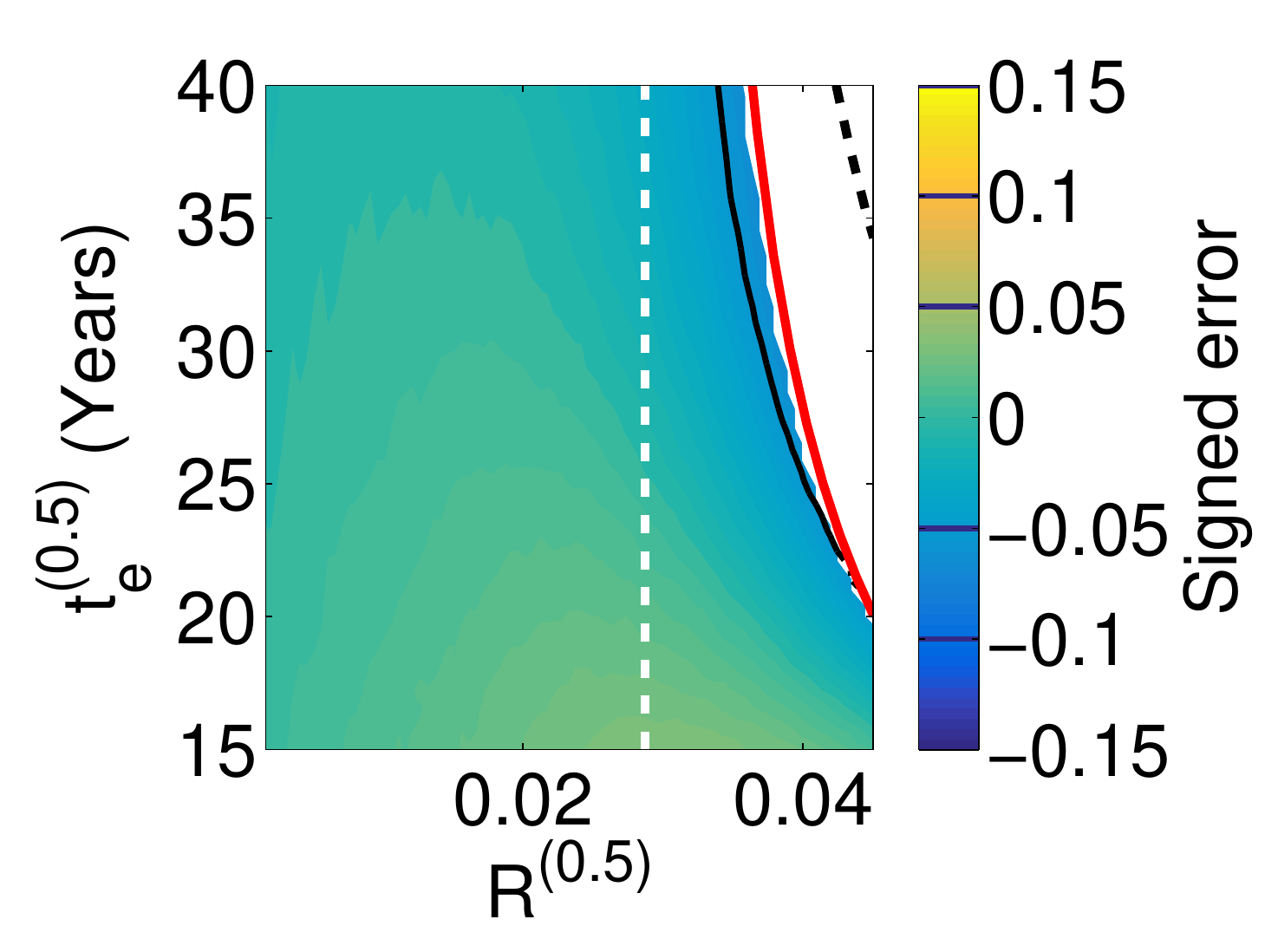}} 
        \hfill
        \subcaptionbox{\label{err slow}}[0.32\linewidth]
                {\includegraphics[width=0.33\linewidth]{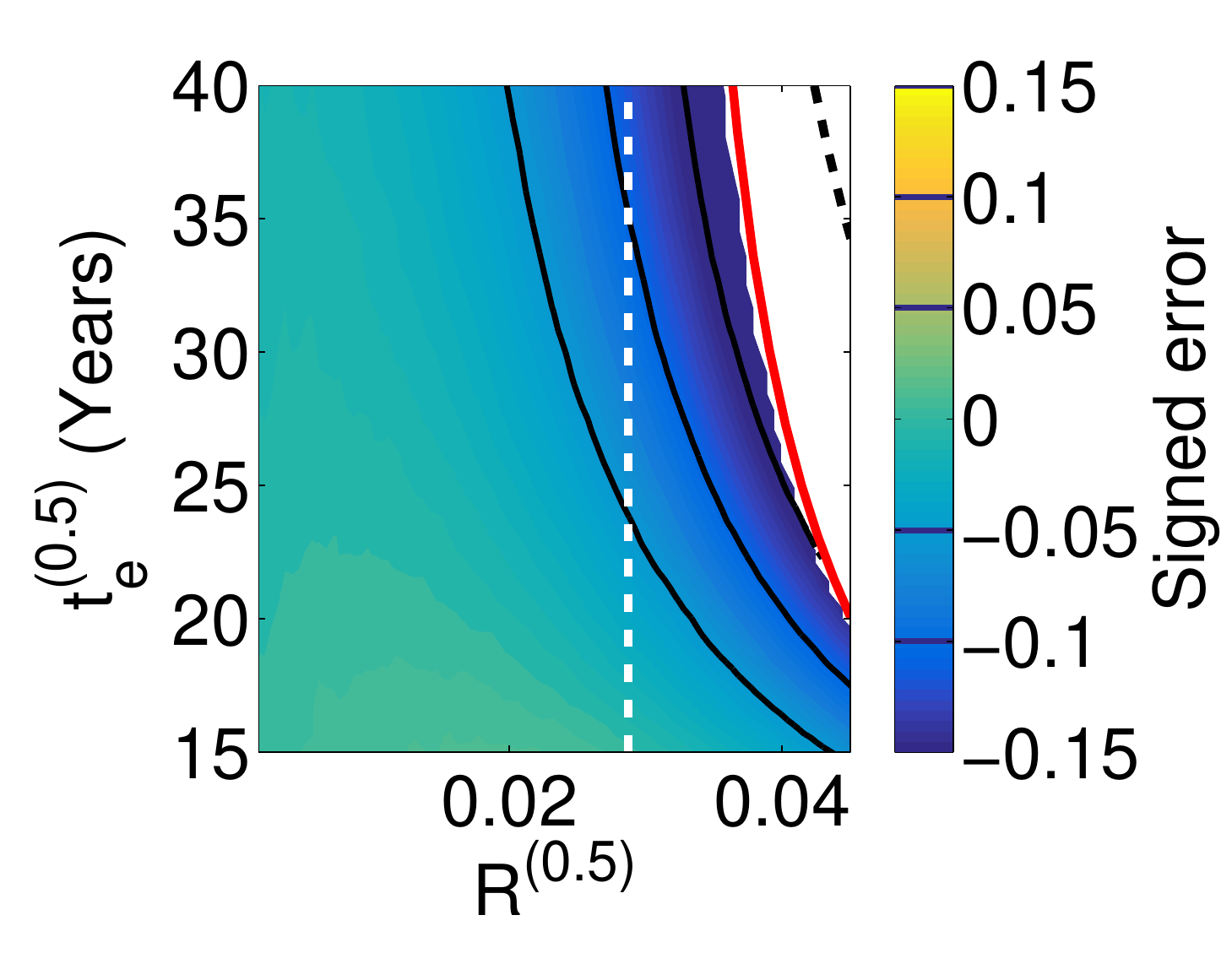}}
                \caption{Comparison between the escape probability as
                  computed from the full two-dimensional model to its
                  projection and approximation \eqref{eq:mode:p}
                  \prrsa{and slow drift approximation (S21) in the
                    Supplementary Material}. (a) Escape probability
                  $P_\mathrm{esc}^\mathrm{2d}$, as computed using
                  \eqref{eq:FP2d} and \eqref{eq:pesc2d}. (b)
                  Difference between mode approximation depicted in
                  Figure~\ref{Eig prob re} of the main article
                  calculated as
                  $P_\mathrm{esc}^\mathrm{1d}-P_\mathrm{esc}^\mathrm{2d}$
                  and panel (a). \prrsa{(c) Difference between slow
                    drift approximation of escape probability,
                    $P_\mathrm{esc}^\mathrm{slow}$, and
                    $P_\mathrm{esc}^\mathrm{2d}$.}  White dashed line
                  indicates the fold bifurcation, black dashed curve
                  the critical parameter values for deterministic
                  tipping and the red solid curve indicates the
                  boundary for validity of the mode approximation.}
  \label{fig:2dapprox}
\end{figure}

The projection of the \prrsa{SDE \eqref{eq:ode+noise}} to a one-dimensional
system gives only accurate predictions for sufficiently small
$\epsilon$, that is, for sufficiently small noise and slow parameter
drift. We compare the predictions from the mode approximation of the
projection onto the atmospheric temperature $T_a$ to the results of
the two-dimensional monsoon model in the main paper, (3.1)--(3.2).  To
evaluate the tipping probability for the two-dimensional monsoon model
we solve the Fokker-Planck equation
\begin{equation}\label{eq:FP2d}
  \begin{aligned}
    \partial_tu=&D_1\partial_{Q_a}^2u+D_2\partial_{T_a}^2u-
    \partial_{Q_a}[f_1(Q_a,T_a,\Asys(t))u]\\
    &-\partial_{T_a}[f_2(Q_a,T_a,\Asys(t))u]
  \end{aligned}
\end{equation}
on the decade timescale with noise variances $D_1=0.01$ and $D_2=3$,
and Dirichlet boundary conditions on the domain
$(Q_a,T_a)\in[-0.04,0.07]\times[295,315]$ (so slightly larger than the
physically realistic ranges). The \prrsa{forcing}s $\Asys(t)$ are chosen as
described in Section 3 and Section 5 of the main article, namely
\begin{equation}
  \Asys(t) = 
  \Asys^{\infty} +  \frac{R^{(0.5)}+\Asys^\mathrm{th}-\Asys^\infty}{\cosh(\prrsa{S}(t_{\mathrm{end}} - 2t))^2}\mbox{}
\label{app:Albedo forcing repeat}
\end{equation}  
with $\Asys^\infty=0.47$ (the approximate present
day value), and $R$ and $\prrsa{S}$ varying such that we exceed the
threshold $\Asys^\mathrm{th}=0.5$ for time $t_e^{(0.5)}\in[15,40]$ years and
amplitude $R^{(0.5)}\in[0,0.045]$. For each simulation, we start from
the eigenvector for the dominant (close to $1$) eigenvalue of the
Fokker-Planck operator on the right-hand side of \eqref{eq:FP2d}. The
total simulation time period is chosen such that, after the transient
time period, $\Asys$ starts from close to its current day value,
namely $0.471$, and returns back to this value again at the end of the
total simulation period. The escape probability is then calculated as
\begin{equation}\label{eq:pesc2d}
  P_\mathrm{esc}^\mathrm{2d}=1-\int_{Q_a,T_a}
  u(Q_a,T_a,t_\mathrm{end})\d Q_a \d T_a,
\end{equation}
at the end of an overall integration period (neglecting the escape during the transient time period). The resulting escape probabilities are shown in
Figure~\ref{True}.
The one-dimensional projection of the monsoon model \eqref{ODE projection} is extracted from the output, temperature $T_a$, using
a second-order approximation of the known equilibrium curve
$(Q_a,T_a,\Asys)$ and the attraction rate toward stable equilibria
nearby. In practice, these quantities may have to be estimated from
observations or model outputs. Then we use a modification of
\prrsa{(S16) from the Supplementary Material}, replacing the parabolic \prrsa{forcing} with the
(rescaled) forcing of the planetary albedo \eqref{app:Albedo forcing repeat}:
\begin{equation}
  \label{eq:mode:nf:comp}
  \dot x=D^{-2/3}a_0^2\kappa(\Asys(D^{-1/3}t)-\Asys^b)+x^2\mbox{.}
\end{equation}
In Eq.\,\eqref{eq:mode:nf:comp} the \prrsa{forcing} $\Asys(t)$
is given by \eqref{app:Albedo forcing repeat}, $\Asys^b\approx0.5287$,
$D=\w_0^T\operatorname{diag}(0.01,3)\w_0=3.04$ with
$\w_0=(-3.50,-0.99)^T$, $\kappa\approx6\cdot10^{-3}$ and
$a_0\approx115.30$ ($a_0^2\kappa$ was estimated from decay rate and
equilibrium curve of the 2d monsoon model.

The difference between the true escape probability (Figure~\ref{True})
and the mode approximation for the one-dimensional projected system
(Figure~\ref{err mode}) is less than $0.05$ in absolute value
everywhere in the
region. 
The main source of error is that, due to the large noise variance, the
system visits parts of the phase space where the quadratic
approximation to the fold and the projection onto a single dimension
are not accurate (the time scale separation between the two dimensions
is only large close to the fold). \prrsa{For comparison,
  Figure~\ref{err slow} shows the slow-drift approximation, which is
  simply accumulating the escape rate for fixed parameter $\Asys$ over
  the interval (see Supplementary Material for precise formulas). We
  observe that the mode approximation has a systematically smaller
  error especially close to the fixed-$\Asys$ bifurcation point (white
  dashed line).}

\section{Monsoon parameters}
\label{app:monsoon vals}

Section \ref{sec:Monsoon} discussed a simplification to the Indian summer monsoon model used by \citet{zickfeld2004modeling}, retaining the key dynamics behind the mechanisms of the monsoon. Table \ref{Table of Parameters} lists all the parameters and their values used in the simplified monsoon model.
\begin{table}[ht]
\centering
\caption{\label{Table of Parameters}Table of parameters used in Indian summer monsoon model}
\def\arraystretch{1.2}
\begin{tabular}{ l c r l}
    Parameter & Value & Unit & Description\\[0.5ex] \hline
     \\[-3ex] 
               $T_{oc}$ & $300$ & $K$ & temperature over Indian Ocean\\    
               $T_0$ & $273.2$ & $K$ & freezing point\\
               $Q_{oc}$ & $0.0190$ & $1$ & humidity over Indian Ocean\\
               $Q_{\mathrm{sat}}$ & $0.0401$ & $1$ & saturated humidity\\
               $\mathcal{L}$ & $2.5\times 10^6$ & $m^2s^{-2}$ & latent heat\\
               $C_E$ & $3.4375\times 10^{-4}$ & $mm\,s^{-1}K^{-1}$ & evaporation factor\\
               $C_P$ & $0.0027$ & $mm\,s^{-1}$ & precipitation factor\\
               $C_{mo}$ & $6.9021\times 10^{-4}$ & $mm\,s^{-1}K^{-1}$ & moisture advection factor ocean\\
               $C_{ml}$ & $1.6213\times 10^{-4}$ & $mm\,s^{-1}K^{-1}$  & moisture advection factor land\\
               $C_{L1}$ & $1.6642$ & $Kg\,s^{-3}K^{-1}$ & outgoing long-wave radiation factor \\
               $C_{L2}$ & $-263.3753$ & $Kg\,s^{-3}$ & outgoing long-wave radiation constant\\
               $F_{\downarrow}^{SL,TA}$ & $443.6250$ & $Kg\,s^{-3}$ & incoming short wave radiation factor\\
               $C_H$ & $0.7136$ & $Kg\,s^{-3}\,K^{-2}$ & heat advection factor\\
               $\theta_{oc}$ & $300.2356$ & $K$ & potential temperature over Indian Ocean\\
               $\Gamma_0$ & $0.0053$ & $K\,m^{-1}$ & atmospheric lapse rate constant\\
               $\Gamma_1$ & $5.5\times 10^{-5}$ & $m^{-1}$ & atmospheric lapse rate linear factor \\
               $\Gamma_2$ & $1000$ & $1$  & atm. lapse rate quadratic factor for $Q_a$\\
               $\Gamma_a$ & $0.0098$ & $K\,m^{-1}$ & rate of decrease for pot. temp. over land\\
               $z_h$ & $5.1564\times 10^3$ & $m$ & high altitude fixed in model\\
               $I_q$ & $2.0636\times 10^3$ & $mm$ & humidity scaling factor \\
               $I_T$ & $1.1958\times 10^9$ & $Kg\,s^{-2}\,K^{-1}$ & temperature scaling factor \\
               $\beta$ & $3.1710\times 10^{-9}$ & decades$\,s^{-1}$ & time scaling factor
    \end{tabular}
\end{table}


\end{appendix}


\bibliography{Monsoon_bib,Rate-induced_paper2}
\bibliographystyle{unsrtnat}
%
%
%
%
%

\end{document}